\documentclass[a4paper]{article}
\usepackage{RR}
\usepackage{hyperref}
\usepackage{graphicx}
\usepackage{amssymb}
\usepackage{amscd}
\usepackage{psfrag}
\newtheorem{remark}{Remark}[section]
\def\begrem{\begin{remark}\rm}
\def\endrem{\null\hfill\blackbox\end{remark}}

\def\RR{{\rm I\hspace{-0.50ex}R} }
\def\be{\begin{equation}}
\def\ee{\end{equation}}
\let\pa\partial
\def \intdb{\int\!\!\!\!\int}

\RRdate{Novembre 2008}

\RRauthor{Nicolas Crouseilles\thanks{INRIA-Nancy-Grand Est, CALVI Project}
  \and 
Thomas Respaud
  \thanks{IRMA Strasbourg et INRIA-Nancy-Grand Est, CALVI Project} 
\and 
Eric Sonnendr\"{u}cker
  \thanks{IRMA Strasbourg et INRIA-Nancy-Grand Est, CALVI Project}
}
\authorhead{Crouseilles, Respaud \& Sonnendr\"{u}cker}
\RRtitle{Une m\'ethode semi-Lagrangienne en avant pour la r\'esolution num\'erique de l'\'equation de Vlasov}
\RRetitle{A Forward semi-Lagrangian Method for the Numerical Solution of the Vlasov Equation}
\titlehead{Forward semi-Lagrangian Vlasov solvers}
\RRresume{Ce document concerne la r\'esolution num\'erique de l'\'equation de Vlasov qui est un 
mod\`ele cin\'etique permettant de d\'ecrire l'\'evolution d'un plasma. Elle est coupl\'ee \`a une \'equation,
par exemple l'\'equation de Poisson, permettant le calcul des champs auto-consistants.
Le mod\`ele coupl\'e est non lin\'eaire. Nous proposons ici une nouvelle m\'ethode semi-Lagrangienne dans laquelle les caract\'eristiques sont int\'egr\'ees en avant, dans le sens du temps, contrairement \`a la m\'ethode semi-Lagrangienne classique qui int\`egre les caract\'eristiques en arri\`ere.
L'autre ingr\'edient de la m\'ethode semi-Lagrangienne est une technique de d\'eposition de l'information port\'ee par les particules sur une grille. Celle-ci est bas\'ee sur un produit tensoriel de splines cubiques de sorte qu'en deux dimensions la formule de d\'eposition implique 16 points de la grille.
Gr\^ace \`a cette nouvelle technique la m\'ethode semi-Lagrangienne devient compl\`etement explicite et peut s'appuyer sur des solveurs num\'eriques d'\'equations diff\'erentielles classiques, rendant ainsi plus simple la mont\'ee en ordre en temps. La m\'ethode est valid\'ee sur plusieurs cas tests repr\'esentatifs des probl\`emes r\'ealistes bas\'es sur ce mod\`ele.
}
\RRabstract{This work deals with the numerical solution of the 
Vlasov equation. This equation gives a kinetic description of the evolution of a plasma, and is coupled with Poisson's equation for the computation of the self-consistent electric field. The coupled model is non linear.
A new semi-Lagrangian method, based on forward integration of the characteristics, is developed. 
The distribution function is updated on an eulerian grid, and the pseudo-particles 
located on the mesh's nodes follow the characteristics of the equation forward for 
one time step, and are deposited on the 16 nearest nodes. This is an explicit way of 
solving the Vlasov equation on a grid of the phase space, which makes it easier to develop 
high order time schemes than the backward method.
}
\RRmotcle{m\'ethode semi-Lagrangienne, Runge-Kutta, simulation de plasmas, \'equation de Vlasov}
\RRkeyword{Semi-Lagrangian method, Runge-Kutta, plasma simulation, Vlasov equation}
\RRprojet{Calvi}
\RRtheme{\THNum } 
 \URLorraine 

\begin{document}
\RRNo{6727}
\makeRR



\section{Introduction}
Understanding the dynamics of charged particles in a plasma is of great importance 
for a large variety of  physical phenomena, such as the confinement 
of strongly magnetized plasmas, or laser-plasma interaction problems 
for example. Thanks to recent 
developments in computational science and in numerical methods, meaningful comparisons 
between experience and numerics are becoming possible. 

An accurate model for the motion of charged particles, is given by the Vlasov equation. It is based on a phase space description so that non-equilibrium dynamics can accurately be investigated. The unknown $f(t, x, v)$ depends 
on the time $t$, the space $x$ and the velocity $v$. The electromagnetic fields 
are computed self-consistently through the Maxwell or Poisson equations, which 
leads to the nonlinear Vlasov-Maxwell or Vlasov-Poisson system.     

The numerical solution of such systems is most of the time
performed using Particle In Cell (PIC) methods, in which 
the plasma is approximated by macro-particles (see \cite{birdsall}). They are 
advanced in time with the electromagnetic fields which are 
computed on a grid. However, despite their capability to treat 
complex problems, PIC methods are inherently noisy, which becomes 
problematic when low density or highly turbulent regions are studied. 
Hence, numerical methods which discretize the Vlasov equation on a 
grid of the phase space can offer a good alternative to 
PIC methods (see \cite{cheng, filbet1, filbet2, sonnen, vecil}).   
The so-called Eulerian methods can deal with strongly nonlinear 
processes without additional complexity, and are well suited for parallel 
computation (see \cite{virginie}). Moreover, semi-Lagrangian methods 
which have first been introduced in meteorology 
(see \cite{staniforth, ZerWoS05, ZerWoS06}), 
try to take advantage of both 
Lagrangian and Eulerian approaches. Indeed, they allow 
a relatively accurate description of the phase space using a fixed mesh 
and avoid traditional step size restriction using the invariance 
of the distribution function along the trajectories. 
Standard semi-Lagrangian methods calculate departure points 
of the characteristics ending at the grid point backward in time; an interpolation step 
enables to update the unknown.

In this work, we consider the numerical resolution 
of the two-dimensional Vlasov equation on a mesh of the phase space 
using a forward semi-Lagrangian numerical scheme. 
In the present method, the characteristics curves are advanced in time and 
a deposition procedure on the phase space grid, similar to the procedure used in PIC methods for the configuration space only, enables to update the distribution function.  
 
One of the main cause for concern with semi-Lagrangian methods, is 
computational costs. With Backward semi-Lagrangian methods (BSL), 
the fields have to be computed iteratively, with Newton fixed point methods, 
or prediction correction algorithms. This is due to an implicit way of solving the 
characteristics (see \cite{sonnen} for details). This strategy makes 
high order resolution quite difficult and expensive. Making the problem 
explicit enables to get rid of iterative methods 
for the characteristics, and to use for example high order Runge Kutta methods 
more easily. This is one of the main advantages of the present Forward semi-Lagrangian 
(FSL) method. Once the new position of the particles computed, a remapping 
(or a deposition) step has to be performed. This issue is achieved using 
cubic spline polynomials which deposit the contribution of the Lagrangian particles 
on the uniform Eulerian mesh. This step is similar to the deposition step 
which occurs in PIC codes but in our case, the deposition is performed in all 
the phase space grid. Similarities can also be found in strategies developed in 
\cite{Cotter, nair, Reich} for meteorology applications. 

In order to take benefit from the advantages of PIC and semi-Lagrangian methods, 
and since the two methods (PIC and FSL) really look like each other, except the 
deposition step, we have also developed a hybrid method, 
where the deposition step is not performed at each time step, but every $T$ time steps. 
During the other time steps, the fields are computed directly 
at the new position of the particles. It shall be noticed that 
however the present method is not a real PIC method, since 
the particle weights are not constant. Indeed, in this method based on 
a description of the unknown using cubic spline polynomials, 
the spline coefficients play the role of the particle weights, 
and are updated at each time step. This kind of hybrid approach 
has been developed recently in a slightly different framework in \cite{parker} 
inspired by \cite{denavit}. 

This paper is organized as follows. In the next section, the two Vlasov equations 
which will be dealt with are presented. Then, we shall introduce the Forward 
semi-Lagrangian (FSL) method, always regarding it comparatively to Backward 
semi-Lagrangian (BSL) methods. Afterward, numerical results for several test 
cases are shown and discussed. Eventually, some specific details are given 
in two appendices, one for the computation of an exact solution to the Landau damping
problem, and the other about the solution of the Poisson equation for the Guiding-Center model.

\section{Models in plasma physics}
\label{model}
In this section, we briefly present two typical reduced models from plasma physics for the description of the time evolution of charged particles.
These two-dimensional models are relevant for more complex 
problems we are interested in and shall be used to validate our new method.

\subsection{Vlasov-Poisson model}
We consider here the classical 1D Vlasov-Poisson model,  
the unknown of which $f=f(t, x, v)$ 
is the electron distribution function. It depends on the space
variable $x\in[0, L]$ where $L>0$ is the size of the domain, 
 the velocity variable $v\in\RR$ and the time 
$t\geq 0$. The Vlasov equation which translates the invariance of the 
distribution function along the characteristics then writes 
\be
\label{vlasov}
\frac{\pa f}{\pa t}+v\partial_x f + E(t, x)\partial_v f= 0,
\ee
with a given initial condition $f(0,x,v)=f_0(x,v)$. 
The self-consistent electric field $E(t, x)$ is computed thanks 
to the distribution function $f$ 
\be
\label{poisson}
\partial_x E(t, x)= \int_{\RR}f(t, x, v)dv - \rho_i,\;\;\;\; \int_0^L E(t,x){\rm d}x=0, 
\ee
where $\rho_i$ denotes the ion density which forms a uniform 
and motionless background in the plasma.

The Vlasov-Poisson model constitutes a nonlinear self-consistent 
system as the electric field determines $f$ with (\ref{vlasov}) 
and is in turn determined by it in (\ref{poisson}). It presents 
several conserved quantities as the total number of particles, 
the $L^p$ norms ($p\geq 1$) defined by $\|f\|_{L^p}=(\intdb |f|^p dx dv)^{1/p}$, 
the momentum and the total energy, as follows:
\begin{eqnarray*}
\frac{d }{d t}\intdb f(t, x, v)dxdv = \frac{d}{d t} \|f(t)\|_{L^p} &=& \frac{d }{d t}\intdb v f(t, x, v)dxdv \nonumber\\ 
&=& \frac{d }{d t}\left[\intdb v^2 f(t, x, v)dxdv + \int E(t, x)^2 dx dv\right] \nonumber\\
&=& 0. 
\end{eqnarray*}
One of the main features of the present work is to develop 
accurate numerical methods which are able to preserve 
these exactly or approximately these conserved quantities for long times. 

\subsection{Guiding-center model}
We are also interested in other kinds of Vlasov equations. 
For instance, in the guiding-center approximation. Charged particles in 
magnetized tokamak plasmas can be modeled by the density $f=f(t, x, y)$ 
in the $2$ dimensional poloidal plane by 
\be
\label{gc}
\frac{\pa f}{\pa t}+E^{\perp}(x, y)\cdot\nabla f = 0,
\ee
coupled self-consistently to Poisson's equation for the electric field 
which derives from a potential $\Phi=\Phi(x, y)$ 
\be
\label{poisson2d}
-\Delta \Phi(t, x, y)= f(t, x, y), \;\;\; E(t, x, y)=-\nabla \Phi(t, x, y). 
\ee
In equation (\ref{gc}), the advection term $E^{\perp} = (E_y, -E_x)$ 
depends on $(x, y)$ and the time-splitting cannot be simply applied like 
in the Vlasov-Poisson case. Hence, this simple model appears to be interesting 
in order to test numerical methods.  

The guiding-center model (\ref{gc})-(\ref{poisson2d}) also 
presents conserved quantities as the total number of particles 
and $L^2$ norm of $f$ (energy) and $E$ (enstrophy) 
\be
\label{l2norm}
\frac{d}{d t}\intdb f(t, x, y)dxdy = \frac{d \|f(t)\|^2_{L^2}}{d t} = \frac{d \|E(t)\|^2_{L^2}}{d t} = 0. 
\ee

\subsection{Characteristic curves}
We can re-write Vlasov equations in a more general context by introducing 
the characteristic curves 
\be
\label{carac}
\frac{d X}{dt} = U(X(t), t). 
\ee
Let us introduce $X(t,x,s)$ as the solution of this dynamical system, 
at time $t$ whose value is $x$ at time $s$. These are called the 
characteristics of the equation. With $X(t)$ a solution of (\ref{carac}), 
we obtain:
\be
\frac{d}{dt}(f(X(t),t)) = \frac{\pa f}{\pa t} + \frac{dX}{dt} \cdot \nabla_X f \\ 
 = \frac{\pa f}{\pa t} + U(X(t),t) \cdot \nabla_X f = 0.
 \ee
which means that $f$ is constant along the characteristics. 
Using these notations, it can be written 
$$
f(X(t;x,s),t) = f(X(s;x,s),s) = f(x,s)
$$
for any times $t$ and $s$, and any phase space coordinate $x$.  
This is the key property used to define semi-Lagrangian methods for the 
solution of a discrete problem.

\section{The forward semi-Lagrangian method}
In this section, we present the different stages of the forward 
semi-Lagrangian method (FSL) and try to emphasize the differences 
with the traditional backward semi-Lagrangian method (BSL).

\subsection{General algorithm}

Let us consider a grid of the studied space (possibly phase-space)
 with $N_{x}$ and $N_{y}$ the number of points in the $x$ direction 
 $[0, L_{x}]$ and in the $y$ direction $[0, L_y]$. We then define 
$$
\Delta x = L_x/N_x, \;\; \Delta y = L_y/N_y, \;\; x_i = i\Delta x, \;\; y_j=j\Delta y, 
$$ 
for $i=0, .., N_x$ and $j=0, .., N_y$. 
One important point of the present method is the definition 
of the approximate distribution functions which are projected on a cubic B-splines basis:
\be
\label{fsl}
f(t, x, y) =  \sum_{k,l} \omega^n_{k,l} S(x-X_1(t;x_k, y_l, t^n))S(v-X_2(t; x_k, y_l, t^n)), 
\ee
where $X(t;x_k, y_l, t^n)=(X_1, X_2)(t, x_k, y_l, t^n)$ 
corresponds to the solution of the characteristics at time $t$ 
(of the two dimensional system (\ref{carac})) 
whose value at time $t^n$ was the grid point $(x_k, y_l)$. 
The cubic B-spline $S$ is defined as follows 
$$
6 S(x)=
\left\{
\begin{array}{ll}
(2-|x|)^3 & \textrm{if} \ 1\le |x|\le 2,\\
4-6x^2+3|x|^3 &\textrm{if}  \ 0\le |x|\le 1,\\
0& \textrm{otherwise.}
\end{array}
\right.
$$
In the expression (\ref{fsl}), the weight $w^n_{k,l}$ 
is associated to the particle located at the grid point 
$(x_k, y_l)$ at time $t^n$; it corresponds to the coefficient of 
the cubic spline determined by the following interpolation conditions
\begin{eqnarray*}
f(t^n, x_i, y_j)&=&\sum_{k,l} \omega^n_{k,l} S(x_i-X_1(t^n;x_k, y_l, t^n))S(y_j-X_2(t^n; x_k, y_l, t^n))\nonumber\\ 
		&=&\sum_{k,l} \omega^n_{k,l} S(x_i-x_k)S(y_j-y_l). 
\end{eqnarray*}
Adding boundary conditions (for example the value of the normal derivative of $f$ at the boundaries, we obtain a set of linear systems in each direction from which the weights
$\omega^n_{k,l}$ can be computed as in 
\cite{sonnen, virginie}.  

We can now express the full algorithm for the forward semi-Lagrangian method
\begin{itemize}
\item Step 0: Initialize $f^0_{i,j}=f_0(x_i,y_j)$ 
\item Step 1: Compute the cubic splines coefficients $\omega^0_{k,l}$ 
such that 
$$
f^0_{i,j} = \sum_{k,l} \omega^0_{k,l} S(x_i-x_k)S(y_j-y_l),  
$$
\item Step 2: Integrate (\ref{carac}) from $t^n$ to $t^{n+1}$, given as initial 
data the grid points $X(t^n)=(x_k, y_l)$ to get $X(t;x_k, y_l, t^n)$ for 
$t\in [t^n, t^{n+1}]$, assuming the advection velocity $U$ is known. We shall explain
in the sequel how it is computed for our typical examples.
\item Step 3: Project on the phase space grid 
using (\ref{fsl}) with $t=t^{n+1}$ to get $f^{n+1}_{i,j}=f^{n+1}(x_i, y_j)$
\item Step 4: Compute the cubic spline coefficients $\omega^{n+1}_{k,l}$ 
such that 
$$
f^{n+1}_{i,j} = \sum_{k,l} \omega^{n+1}_{k,l} S(x_i-x_k) S(y_j-y_l).   
$$
\item Go to Step 2 for the next time step.
\end{itemize}

\subsection{{FSL: An explicit solution of the characteristics}}
For BSL, especially for the solution of the characteristics, 
it is possible to choose algorithms based on two time steps with 
field estimations at intermediate times. 
Generally, you have to use a fixed-point algorithm, a Newton-Raphson 
method (see \cite{sonnen}), a prediction correction one or also 
Taylor expansions (see \cite{virginie}) in order 
to find the foot of the characteristics. This step of the global 
algorithm  costs a lot (see \cite{sonnen}). 
It is no longer needed in FSL, where the starting point of the 
characteristics is known so that traditional methods to solve ODEs, 
like Runge-Kutta algorithms can be incorporated to achieve high order 
accuracy in time.  Let us show the details of this explicit solution 
of the characteristics, in Vlasov-Poisson and Guiding-Center models.

\vspace{5.mm}

In both cases, the dynamical system (\ref{carac}) has 
to be solved. With FSL, $X(t^n)$, $U(X(t^n),t^n)$ are known. 
You can choose your favorite way of solving this system on 
each time step, since the initial conditions are explicit. 
This leads to the knowledge of $X(t^{n+1})$ and $U(X(t^{n+1}),t^{n+1})$ 
so that Step $2$ of the previous global algorithm is completed.

As examples of forward solvers for the characteristics curves, 
the second-order Verlet algorithm, Runge-Kutta $2$ 
and Runge-Kutta $4$ will be proposed for Vlasov-Poisson, and, as Verlet cannot be applied, only Runge-Kutta $2$ and $4$ will be used for the Guiding-Center model.

For Vlasov-Poisson, we denote by 
$X(t^n)=(X_1(t^n), X_2(t^n))=(x^n, v^n)$ the mesh of 
the phase space, and $U(X(t^n), t^n) = (v^n, E(x^n, t^n))$ the advection velocity.  
The Verlet algorithm can be written  
\begin{itemize}
\item Step $1$: $\forall k,l$, $v_{k,l}^{n+\frac{1}{2}} -v_l^n \,=\,  \frac{\Delta t}{2} \,  E(x_k^n, t^n)$,
\item Step $2$: $\forall k,l$, $x_{k,l}^{n+1} -x_k^n \,=\, \Delta t  \, v_{k,l}^{n+1/2}$,   
\item Step $3$: compute the electric field at time $t^{n+1}$
\begin{itemize}
\item   deposition of the particles $x_{k,l}^{n+1}$ on the spatial grid $x_i$ for the density $\rho$:
$\rho(x_i, t^{n+1})=\sum_{k,l} \omega^n_{k,l} S(x_i-x^{n+1}_{k,l})$, like in a PIC method.
\item solve the Poisson equation on the grid $x_i$: $E(x_i, t^{n+1})$.   
\end{itemize}
\item Step $4$: $\forall k,l$, $v_{k,l}^{n+1} -v_{k,l}^{n+\frac{1}{2}} \,=\, \frac{\Delta t}{2}  \, E(x_{k,l}^{n+1}, t^{n+1})$.   
\end{itemize}

A second or fourth order Runge-Kutta algorithms can also be used to solve 
the characteristics curves of the Vlasov-Poisson system forward in time. 
The fourth order Runge-Kutta algorithm needs to compute intermediate 
values in time of the density and the electric field. Let us detail the algorithm 
omitting the indices $k,l$ for the sake of simplicity 
\begin{itemize}
\item Step $1$: $k_1 = (v^n, E(x^n, t^n)) = (k_1(1), k_1(2))$,   
\item Step $2$: compute the electric field at intermediate time $t_1$: 
\begin{itemize}
\item   deposition of the particles on the spatial grid $x_i$ for the density $\rho$: 
$\rho(x_i, t_1)=\sum_{k,l} \omega^n_{k,l} S[x_i-(x_k^{n}+\Delta t/2 \; k_1(1))]$. 
\item   solve the Poisson equation on the grid $x_i$: $E(x_i, t_1)$.   
\end{itemize}
\item Step $3$: compute $k_2=(v^n+\frac{\Delta t}{2} k_1(2), E(x^n+\frac{\Delta t}{2}k_1(1), t_1)$
\item Step $4$: compute the electric field at intermediate time $t_2$: 
\begin{itemize}
\item   deposition of the particles on the spatial grid $x_i$ for the density $\rho$:
$\rho(x_i, t_2)=\sum_{k,l} \omega^n_{k,l} S[x_i-(x_k^{n}+\Delta t/2 \; k_2(1))]$. 
\item   solve the Poisson equation on the grid $x_i$: $E(x_i, t_2)$.   
\end{itemize}
\item Step $5$: compute  $k_3=(v^n+\frac{\Delta t}{2} k_2(2), E(x^n+\frac{\Delta t}{2}k_2(1), t_2)$ 
\item Step $6$: compute the electric field at intermediate time $t_3$:
\begin{itemize}
\item   deposition of the particles on the spatial grid $x_i$ for the density $\rho$: 
$\rho(x_i, t_3)=\sum_{k,l} \omega^n_{k,l} S[x_i-(x_k^{n}+\Delta t \; k_3(1))]$. 
\item   solve the Poisson equation on the grid $x_i$: $E(x_i, t_3)$.   
\end{itemize}
\item Step $7$: compute $k_4=(v^n+\Delta t \; k_3(2), E(x^n+\Delta t \; k_3(1), t_3)$  
\item Step $8$: $X^{n+1} - X^n = \frac{\Delta t}{6}\left[k_1+2k_2+2k_3+k_4\right]$
\end{itemize}

In both Verlet and Runge-Kutta algorithms, the value of $E$ at intermediate 
time steps is needed (step $3$ for Verlet and steps $3$, $5$ and $7$ 
for Runge-Kutta $4$). This is achieved as in PIC algorithms by advancing the 
particles (which coincide at time $t^n$ with the mesh in this method) up to the required 
intermediate time. Using a deposition step, the density is computed thanks 
to cubic splines of coefficients $w^n_i$ on the mesh 
at the right time, and thus the electric field can also be computed 
at the same time thanks to the Poisson equation. 
Using an interpolation operator, the electric field is then evaluated at the 
required location (in steps $3, 5$ and $7$). 
Let us remark that this step involves a high order interpolation operator 
(cubic spline for example) which has been proved in our experiments 
to be more accurate than a linear interpolation (see section $4$).  

\vspace{5.mm}

For the Guiding-Center equation, the explicit Euler 
method, and also Runge-Kutta type methods (of order $2$, $3$ and $4$)
have been implemented. 
There is no technical difficulty with computing high order methods. 
This is one of the general interests of forward methods. The time algorithm 
for solving the characteristics at the fourth order is similar to 
those presented in the Vlasov-Poisson case. However, 
there is a additional difficulty in the deposition step which enables 
to evaluate the density at intermediate time steps; indeed, the deposition is 
two-dimensional since the unknown does not depend on the velocity variable in this case. 

Let us summarize the main steps of the 
second order Runge-Kutta method applied to the guiding 
center model of variables $X^{n} = (x^n, y^n)$ and of advection field 
$U(X^n, t^n)=E^{\perp}(X^n, t^n)$
\begin{itemize}
\item Step $1$: $\tilde{X}^{n+1} - X^n = \Delta t E^{\perp}(X^n, t^n)$  
\item Step $2$: Compute the electric field at time $t^{n+1}$ 
\begin{itemize}
\item   two-dimensional deposition of the particles on the spatial grid $(x_j, y_i)$ for the density $\rho$:
$\rho(x_j, y_i, t^{n+1})=\sum_{k} \omega^n_k S[x_j-x_{k,l}^{n+1}]S[y_i-y_{k,l}^{n+1}]$
\item   solve the two-dimensional Poisson equation on the grid $x_j$: $E(x_j, y_i, t_{n+1})$.
\end{itemize}
\item Step $3$: $X^{n+1} - X^n = \frac{\Delta t}{2}\left[E^{\perp}(X^n, t^n) + E^{\perp}(\tilde{X}^{n+1}, t^{n+1})\right]$      
\end{itemize}

Here, the numerical solution of the two-dimensional Poisson's equation 
is based on Fourier transform coupled with finite difference method. 
See details in Appendix II.

%

\subsection{{FSL - BSL Cubic Spline Interpolation}}
We are going to compare how spline coefficients are computed 
recurrently, for one dimensional transport problems, 
for the sake of simplicity.

\paragraph{FSL: deposition principle}

On our mesh, the grid points $x_i=i\Delta x, i=0, .., N_x$ 
at a time $n$ can be regarded as particles. 
We have a distribution function which is projected onto a cubic splines basis. 
Thus, we know $f(t^n, x), \; \forall x$, then the particles move forward, 
and we have to compute $f(t^{n+1}, x_i), \; i=0, ..., N_x$, reminding 
that $f$ is constant along the characteristics, and that the particles 
follow characteristics between $t^n$ and $t^{n+1}$. 

In fact, to each mesh point $x_i$, a spline coefficient $ \omega_k$ is linked. 
The thing to understand, is that these coefficients are transported up to 
the deposition phase. The key is then to compute them recurrently as follows: 
\begin{itemize}
\item Deposition step 
\begin{eqnarray*}
f^{n+1}(x_i) &=&  \sum_k \omega^n_k S(x_i-X(t^{n+1};x_k, t^n)) \nonumber\\
	     &=& \sum_{k/X(t^{n+1};x_k, t^n) \in [x_{i-1},x_{i+2}]}  \omega_k^n S(x_i-X(t^{n+1};x_k, t^n)), 
\end{eqnarray*}
\item Update of the splines coefficients $\omega_k^{n+1}$ using the 
interpolation conditions
$$
f^{n+1}(x_i)= \sum_{k=i-1}^{i+2}  \omega_k^{n+1} S(x_i-x_k), 
$$
\end{itemize}
The number of points which actually take part in the new value of 
$f^{n+1}(x_i)$ (here $4$) is directly linked with the 
spline degree you choose. 
A p-Spline for example has a $(p+1)$ points support. 

In a 1D way of regarding the problem, you can easily prove 
the mass conservation:
\begin{eqnarray*}
m^{n+1} &=&\Delta x\sum_{i} f^{n+1}(x_i)\nonumber\\  
        &=&\Delta x \sum_i \sum_k \omega^{n}_k S(x_i-X(t^{n+1}; x_k, t^n)) \nonumber\\  
    	&=& \Delta x\sum_k \omega^n_k = \Delta x\sum_{i} f^{n}(x_i)= m^n. 
\end{eqnarray*}
Merely with the spline property of unit partition $\sum_i S(x-x_i) = 1$ for all $x$.


\paragraph{BSL: interpolation principle}
Let us introduce some notations. The foot of the characteristics  
$X(t^n, x_i, t^{n+1})$ belongs to the interval $[x_l,x_{l+1}[$. 
Then the reconstructed distribution function can be written 
\begin{itemize}
\item Interpolation step 
\begin{eqnarray*}
f^{n+1}(x_i) &=& f^{n}(X(t^n;x_i,t^{n+1})) \nonumber\\ 
             &=& \sum_{k = l - 1}^{l+2}  \omega_k^n S(X(t^n;x_i,t^{n+1})-x_k)
\end{eqnarray*}
\item Update of the splines coefficients $\omega^{n+1}_k$ using the 
interpolation conditions
$$
f^{n+1}(x_i)= \sum_{k=i-1}^{i+2}  \omega_k^{n+1} S(x_i-x_k)
$$
\end{itemize}
Here, we denoted by $X(t^n;x_i,t^{n+1})$ the foot of the characteristic 
coming from $x_i$. The reader is referred to \cite{sonnen, virginie} 
for more details on BSL interpolation. 

In both cases, a linear system has to be solved, of equivalent 
complexity, so our method is as efficient at this level as BSL is.

\subsection{{Basic differences between FSL and BSL}}
Let us now explain the basic differences between forward 
and backward semi-Lagrangian methods. In both cases, 
a finite set of mesh points $(x_m)_{m=1..N}$ is used, and the 
values of the function $f$ at the mesh points at a given time step $t^n$
are considered. The aim is to find the new values of $f$ on the 
grid at the next time step $t^{n+1}$.

\paragraph{BSL}
For BSL, in order to find the $(n+1)$-th value of $f$ at $x_m$, we 
follow the characteristic curve which goes through $x_m$, 
backward in time, until time $t^n$. The arrival point 
will be called the foot of the characteristics and does 
not necessarily coincide with a mesh point. Hence, we use 
any interpolation technique to compute $f$ at this point, knowing
all the values of the mesh at this time. 
This leads to the new value of $f(x_m)$. 
Let us summarize:
\begin{itemize}
\item find the foot of the characteristics $X(t^{n})$ knowing $X(t^{n+1})=x_m$ (mesh point) 
\item interpolate using the grid function which is known at time $t^n$.
\end{itemize}

\paragraph{FSL}
For FSL, the principle is quite different. The 
characteristics beginning at time $t^n$ on the grid points are followed, 
during one time step, and the end of the characteristics 
({\it i.e.} at time $t^{n+1}$) is found. At this moment, the 
known value is deposited to the nearest grid points (depending on the chosen method). 
This deposition step is also performed in PIC codes on the spatial grid only, in order to get the sources for the computation of
the electromagnetic field. 
Once every grid points has been followed, the new value of $f$ 
is obtained by summing all contributions.  The FSL method can be 
summarized as follows 
\begin{itemize}
\item find the end of the characteristics $X(t^{n+1})$ 
leaving from $X(t^{n})=x_m$ (mesh point) 
\item deposit on the grid and compute the new particle weights.
\end{itemize}

\begin{figure}
\includegraphics[width=5cm,height=5.cm]{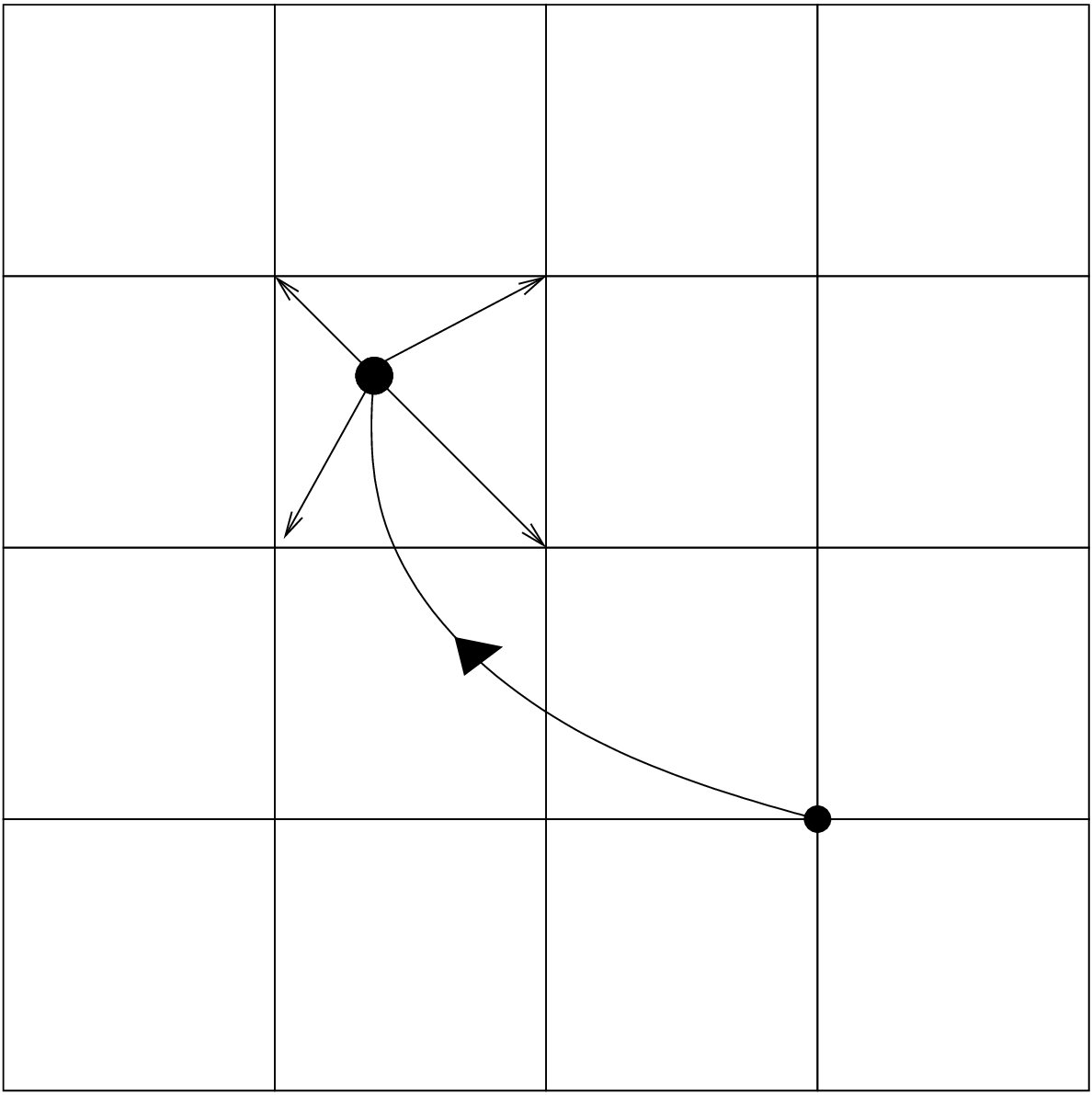}\hfill
\includegraphics[width=5cm,height=5.cm]{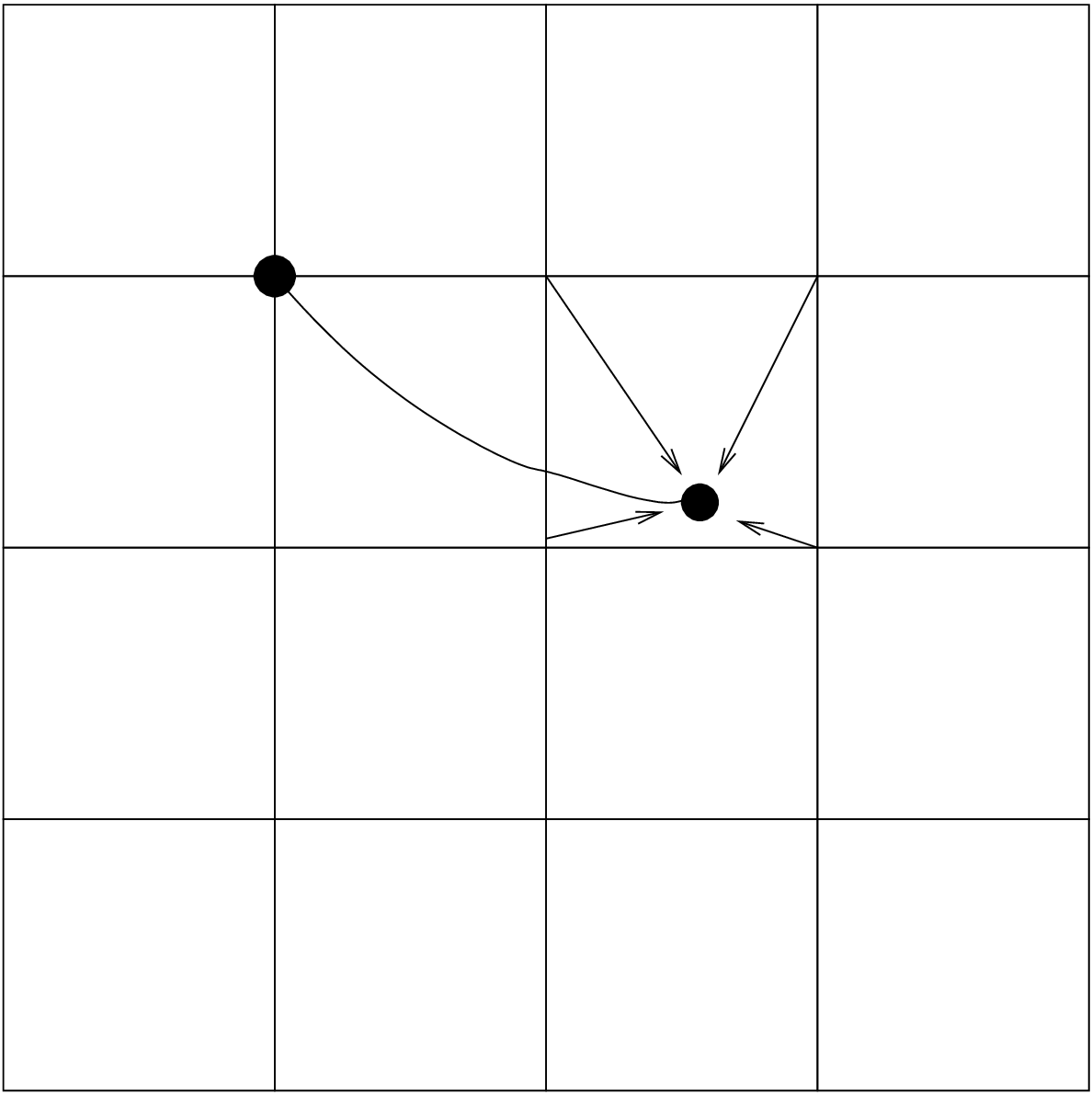}
\caption{Principle of FSL (left) and BSL (right) for linear splines.}
\end{figure}


\section{Numerical results}
This section is devoted to the numerical implementation of 
the forward semi-Lagrangian method. In particular, comparisons 
with the backward semi-Lagrangian method will be performed 
to validate the new approach.

\subsection{Hill's equation}
In order to check that high orders are really reached, a particularly easy model can be used, in which there are no self-consistent fields. This leads to a 1D model with an external  force field written $-a(t)x$, where $a$ is a given periodical function. The Vlasov equation becomes:
\be
\label{vlasov_hill}
\frac{\pa f}{\pa t}+v\partial_x f - a(t)x \partial_v f= 0,
\ee
The solution of this equation is seen through its characteristics, solutions of
\be
\label{char_hill}
\frac{dX}{dt}= V, \;\;\;\;\;\;  \frac{dV}{dt}= -a(t)X
\ee
thus, X is solution of Hill's equation:
\be
\label{hill}
\frac{d^2X}{dt^2} + a(t)X =0
\ee
Let's note that this equation can be written in a general way $\frac{du}{dt}=A(t)u$, 
where $A$ is a matrix valued periodic function. Since this is a 2D linear system, 
its solution is a 2D vector space and it is sufficient to find two independent solutions.


Let $\omega$, $\psi$ $\in C^2(\mathbb{R^+},\mathbb{R})$, with $\omega(t)>0 \;\; \forall t \in \mathbb{R^+}$, so that $\omega$ is solution of the differential equation 

\be
\label{omega}
\frac{d^2\omega}{dt^2} + a(t)\omega - \frac{1}{\omega^3}=0 \;\;\; \frac{d\psi}{dt}=\frac{1}{\omega^2}
\ee   

So $u(t)=\omega(t)e^{i\psi(t)}$ and $v(t)=\omega(t)e^{-i\psi(t)}$ are two independent solutions of Hill's equation (see \cite{magnus} for more details) which can be determined numerically.

\begin{figure}
\psfrag{f}{$x^2 \;$}
\psfrag{g}{$x^4 \;$}
\psfrag{number points}{$\Delta t$}
\begin{tabular}{cc}
\includegraphics[width=7cm,height=7.cm,angle=-0]{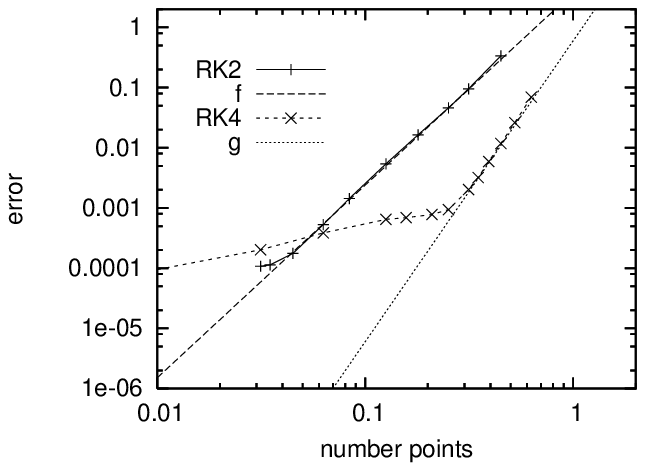} &     
\includegraphics[width=7cm,height=7.cm,angle=-0]{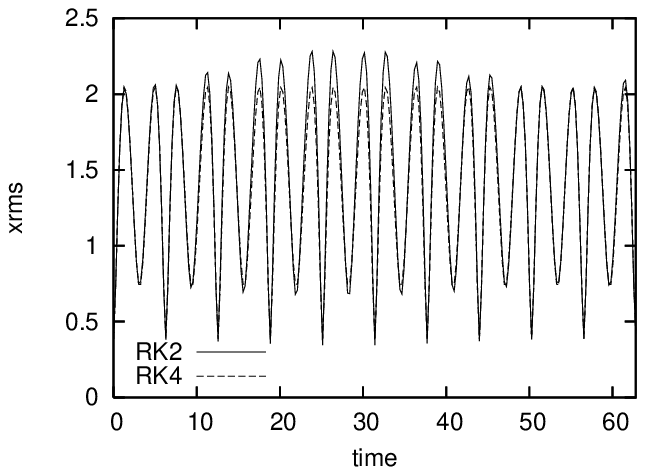}
\end{tabular}
\caption{Error as a function of $\Delta t$ for RK2 and RK4 (left) and $x_{rms}$ as a function of time, for $\Delta t= 2\pi/25$, RK2 and RK4 (right).}
\label{ordre}
\end{figure}

For this test case, the initial distribution function will be:
$$
f_0(x,v)=e^{-\frac{x^2}{2\omega^2}-\frac{\omega^2v^2}{2}}, \forall (x, v)\in [-12,12]^2. 
$$
The associated solution $f(x,v,t)$  will depend only on $A\omega(t)$. In particular, $f$ will have the same periodicity as $a$ and $\omega$. This is what will be used for testing the code.
For different orders (2 and 4), and different $\Delta t$, $x_{rms}(t)=\sqrt{\int{x^2 f(x,v,t) dxdv}}$ will be displayed on Fig \ref{ordre}.  This function should be periodic, and thus should reach the same test value $x_{rms}(0)$ at each period. The error will be measured 
between the ten first computed values and the exact one, for the ten first periods. Then these errors will be summed, so that a $L^1$ norm of the error is dealt with: 
$$
err=\sum_{k=0}^{k=10} e_k, \;\; \mbox{ with } \;\; e_k=|x_{rms}(2k\pi) - x_{rms}(0)|. 
$$
The order of the method is checked in Figure \ref{ordre}. 
Note that $N_x=N_v=1024$ to make sure that 
convergence is achieved for the interpolation step. 
The expected order is achieved for a certain $\Delta t$ interval. 
If $\Delta t$ becomes too small, a kind of saturation happens. 
This is due to the term in $\frac{h^{m+1}}{\Delta t}$ 
(where $h=\Delta x=\Delta v)$ in the theoretical estimation 
of the error for backward methods (\cite{mehrenberger}), which 
becomes too high and prevents us from keeping the correct order. 
A forthcoming paper will try to do the same kind of error estimation for the forward method.

\subsection{Vlasov-Poisson case}
\paragraph{Landau damping}
The initial condition associated to the scaled Vlasov-Poisson equation 
(\ref{vlasov})-(\ref{poisson}) has the following form
\be
\label{initlandau}
f_0(x, v)=\frac{1}{\sqrt{2\pi}}\exp(-v^2/2)(1+\alpha \cos(kx)), \;\;\; (x, v)\in [0, 2\pi/k]\times \RR, 
\ee
where $k=0.5$ is the wave number and $\alpha=0.001$ is the 
amplitude of the perturbation,  
so that linear regimes are considered here. 
A cartesian mesh is used to represent the phase space with a computational domain 
$[0,2\pi/k]\times [-v_{max},v_{max}]$, $v_{max}=6$. The number of mesh points in 
the spatial and velocity directions is designated by $N_x=64$ and $N_v=64$ respectively. 
Finally, the time step is equal to $\Delta t=0.1$ and 
the Verlet algorithm is used to compute the characteristics. 

In this context, it is possible to find an exact value of the dominant mode
solution of the linearized Vlasov-Poisson equation (see Appendix I for some details).
The exact electric field corresponding to the dominant mode reads
$$
E(x,t) =  4 \alpha \times 0.3677 e^{-0.1533 t} \sin(0.5x) \cos(1.4156 t - 0.5326245). 
$$

\begin{figure}
\psfrag{E1}{}
\psfrag{E2}{}
\psfrag{pente}{$\gamma=-0.1533$}
\psfrag{time}{$t \omega_p^{-1}$}
\begin{tabular}{cc}
\includegraphics[width=7cm,height=7.cm,angle=-0]{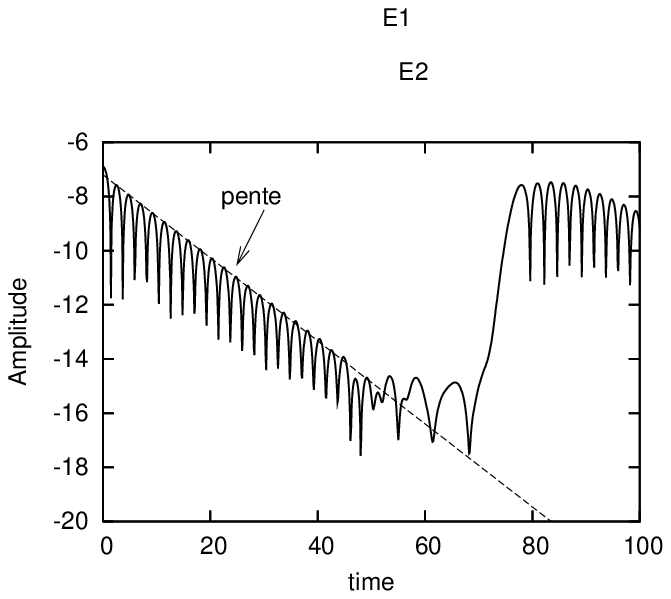} &     
\includegraphics[width=7cm,height=7.cm,angle=-0]{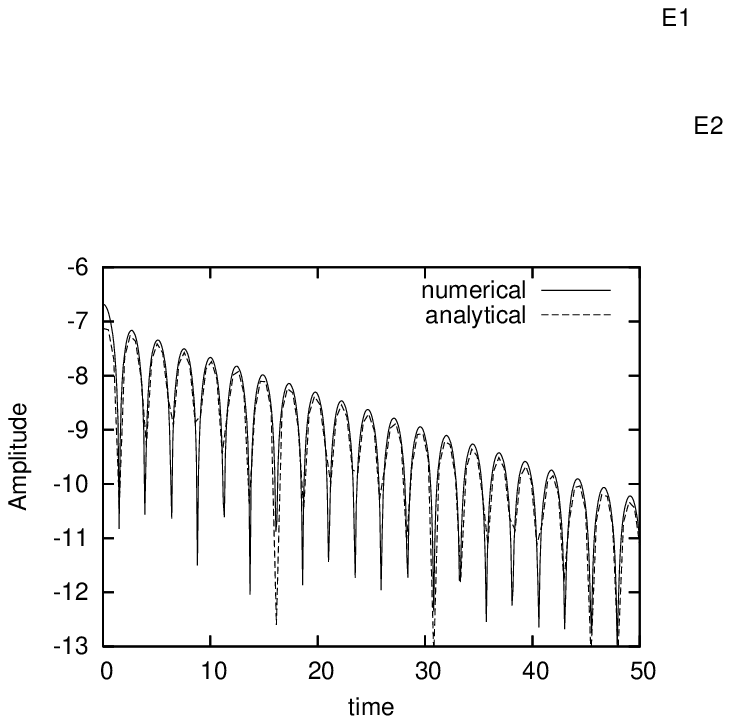}
\end{tabular}
\caption{Linear Landau damping for $k=0.5$ (left) and $k=0.4$ (right)}
\label{landau}
\end{figure}

On Fig. \ref{landau}, the analytical solution of the $L^2$ norm of the  electric field 
and the implemented one are plotted. It can be observed that the 
two curves are very close to    
each other. In particular, the damping rate and the frequency of the wave are 
well recovered ($\gamma=-0.1533$ and $\omega=1.4156$) by the method. Similar 
precision is achieved for different values of $k$ leading to different value of 
the damping rate and of the frequency (see Fig. \ref{landau}). 

The recurrence effect that occurs with the present 
velocity discretization on a uniform grid, at $T_R \approx 80 \, \omega_p^{-1}$ 
can also be remarked. This value is in good agreement with the 
theoretical recurrence time which can be 
predicted in the free-streaming case (see \cite{nakamura}) 
$T_R = \frac{2\pi}{k\Delta v}$.  

\begin{figure}
\psfrag{E1}{}
\psfrag{nb=1}{}
\psfrag{pente}{$\gamma=-0.1533$}
\psfrag{time}{$t \omega_p^{-1}$}
\begin{tabular}{cc}
\includegraphics[width=7cm,height=7.cm,angle=-0]{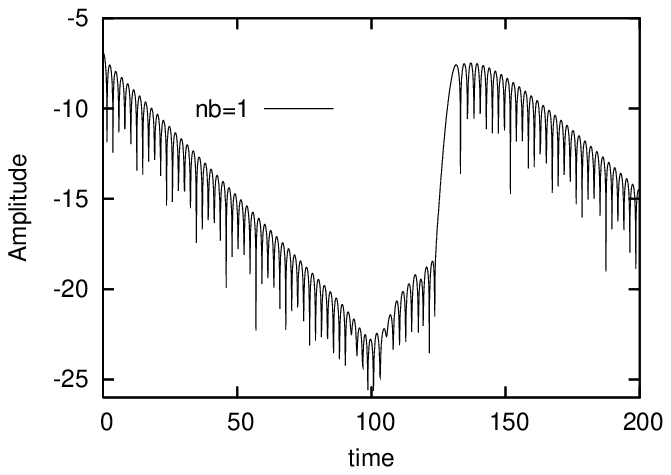} &     
\includegraphics[width=7cm,height=7.cm,angle=-0]{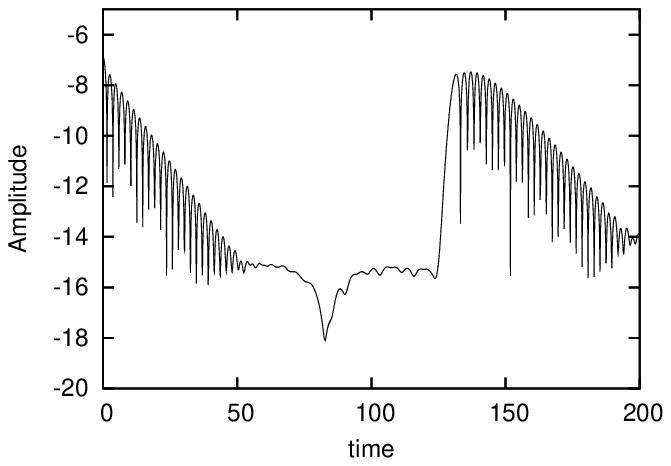}
\\
\includegraphics[width=7cm,height=7.cm,angle=-0]{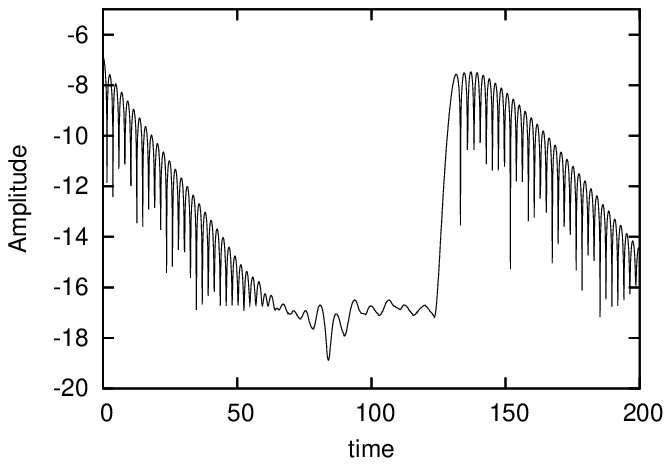} &     
\includegraphics[width=7cm,height=7.cm,angle=-0]{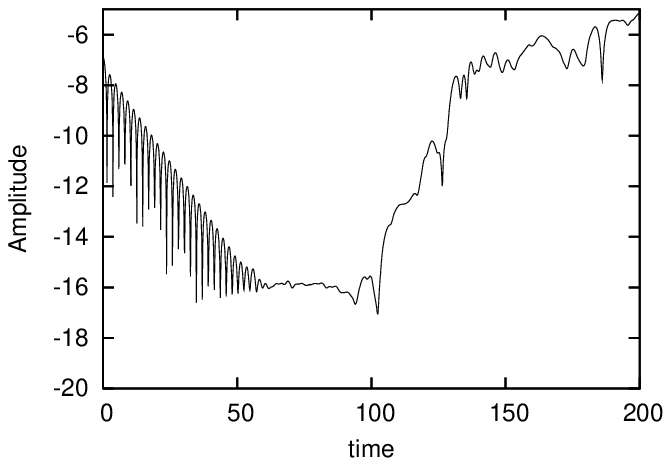}
\end{tabular}
\caption{Linear Landau damping for $k=0.5$ and for different number of $T$: from up to down and left to right: $T=1, T=2, T=16, T=256$.}
\label{landau_deposit}
\end{figure}

This test case has also been solved with the ``hybrid'' method in which the deposition step is 
only performed every $T$ time steps. In all other steps, the remapping (or deposition) step 
is not performed, therefore, it can be linked with a PIC method. As it was already
 said, it is not really a PIC method since the spline coefficients are different on the 
phase space grid and are updated at each remapping step, whereas in classical PIC methods, 
these coefficients (called {\it weights}) are constant equal to 
$\frac{n_0}{N_{part}}$ where $N_{part}$ is the number of particles. 
On Fig. \ref{landau_deposit}, the electric field is plotted again, for 
different values of $T$, and $\Delta t = 0.1$, with $N_x=N_v=128$ points. As expected, 
the method works well, even for large values of $T$. 
Only a kind of 
saturation can be observed, and it can be seen that values smaller than $2^{-18}$ are not 
well treated, because of the lack of accuracy of the hybrid method. Nevertheless, the results 
are convincing: the computation gets faster as $T$ gets larger, 
and a good accuracy is still reached.

\paragraph{Two stream instability}
This test case simulates two beams with opposit 
velocities that encounter (see \cite{filbet2, nakamura}). 
The corresponding initial condition can be given by
$$
f_0(x, v) = {\cal M}(v)v^2[1-\alpha\cos(k x)], 
$$ 
with $k=0.5$ and $\alpha=0.05$. The computational domain is 
$[0,2\pi/k]\times [-9,9]$ which is sampled by $N_x=N_v=128$ points. 
The Verlet algorithm is used to compute the characteristics with $\Delta t=0.5$. 

We are interested in the following diagnostics: the  first three modes 
of the electric field, the electric energy $1/2 \|E(t)\|^2_{L^2}$ 
and the time evolution of the phase space distribution function.  

\begin{figure}
\psfrag{time}{$t \omega_p^{-1}$}
\psfrag{Electric energy}{Amplitude}
\begin{tabular}{cc}
\includegraphics[width=7cm,height=7.cm,angle=-0]{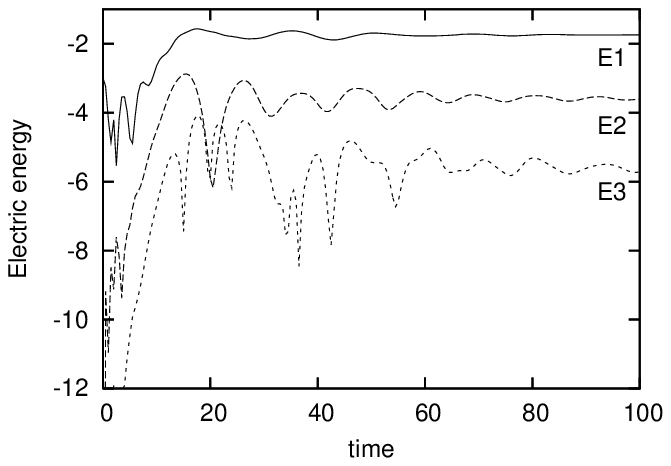} &
\psfrag{Amplitude}{Electric energy}
\includegraphics[width=7cm,height=7.cm,angle=-0]{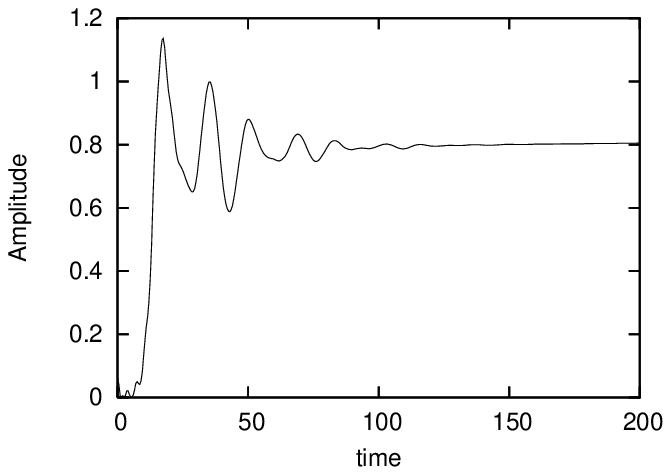}
\end{tabular}
\caption{Two stream instability: Time evolution of the three first modes of the electric field (left) and of the electric energy (right). }
\label{tsi}
\end{figure}

On Fig. \ref{tsi}, we plot the time history of the first three Fourier modes 
of the electric field: $|E_1|, |E_2|, |E_3|$ denotes the amplitudes of 
$\hat{E}(k=0.5)$ , $\hat{E}(k=1)$ and $\hat{E}(k=1.5)$ respectively. 
We observe that after an initial phase, the first mode exponentially increases 
to reach its maximum at $T\approx 18 \, \, \omega_p^{-1}$. 
After this phase and until the end of the simulation, a periodic behavior is 
observed which translates the oscillation of the trapped particles in the electric 
field; in particular, a vortex rotates with a period of 
about $18 \,\omega_p^{-1}$. 
The other modes $|E_2|$ and  $|E_3|$ also grow exponentially and oscillate 
after the saturation. However, their amplitude remains inferior to that of the first mode. 
Similar observations can be performed for the electric energy which reaches  
its maximum at $T\approx 18 \, \omega_p^{-1}$ after an important and fast 
increase (from $t=8$ to $t=18\,\omega_p^{-1}$). 

\begin{figure}
\psfrag{time}{$t \omega_p^{-1}$}
\psfrag{Electric energy}{Amplitude}
\psfrag{nb}{T}
\begin{tabular}{cc}
\includegraphics[width=7cm,height=7.cm,angle=-0]{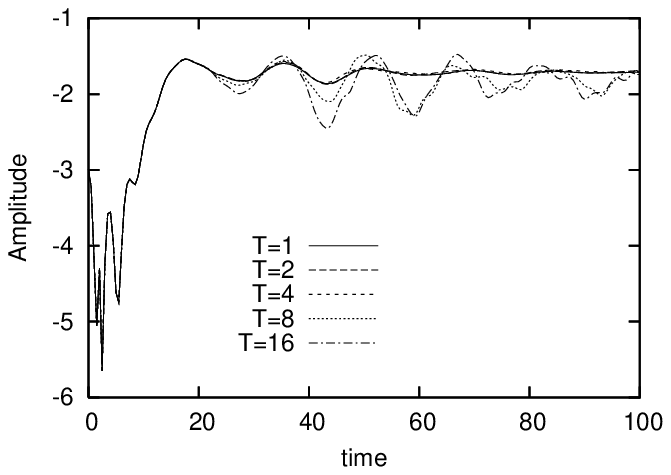} &
\psfrag{Amplitude}{Electric energy}
\includegraphics[width=7cm,height=7.cm,angle=-0]{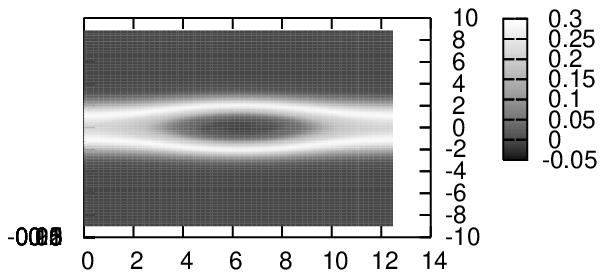}
\\
\includegraphics[width=7cm,height=7.cm,angle=-0]{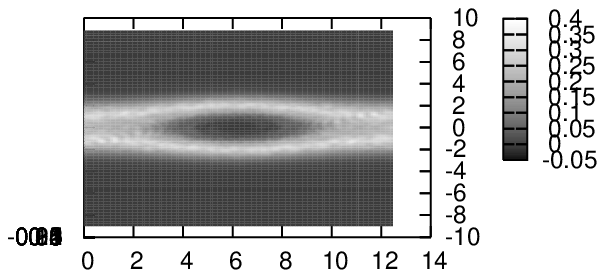} &
\includegraphics[width=7cm,height=7.cm,angle=-0]{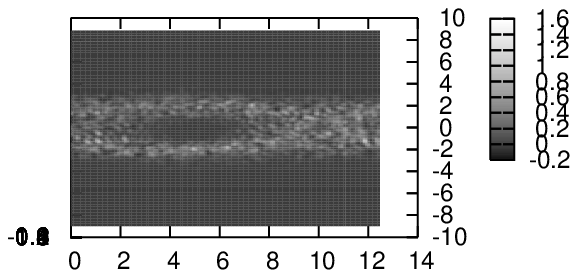}
\end{tabular}
\caption{Two stream instability: Time evolution of the first mode of the electric field (up and left) and distribution function at time $t=100\, \omega_p^{-1}$ for $T=1, 4, 8$.}
\label{tsi_deposit}
\end{figure}

This test case was also solved with the hybrid method to test the capability 
in the nonlinear regime. 
On Fig. \ref{tsi_deposit}, the first Fourier mode of the electric field is displayed for 
different $T$, with $\Delta t = 0.5$, and $128$ points in each direction. It can be observed 
that during the first phase, which is a linear one, even for big $T$, the results are quite 
accurate for all values of $T$. 
The hybrid method seems to have more difficulty after this linear phase. 
Numerical noise, one of the main drawbacks of PIC methods can be observed as $T$ gets higher. 
The phenomenon can be understood looking at the distribution function. On Fig. \ref{tsi_deposit} 
the noise clearly appears on $f$. Noisy values quickly reach 
high values which prevent the method from being accurate enough. They are more important in our 
hybrid method than in classical PIC ones, because of the deposition step, where the particles 
weight play their role. Indeed, if particles with very different weights are located at the 
same place, the deposition does not take into account properly the particles of low weights 
compared to those of heavy weights. On Fig. \ref{tsi_deposit}, it can 
be seen that the vortex which appears at the middle of the distribution and should 
stay there slowly leaves out of the domain. This can be explained as a kind of diffusion. 
Nevertheless, if $T$ remains very little $(2-4)$, the results are really good. It can also be noticed 
that $\Delta t$ plays an important role, actually, when a smaller $\Delta t$ is chosen, the results 
remain good for bigger $T$. As an example, if $\Delta t = 0.1$, results 
remain acceptable until $T=16$. 


\paragraph{Bump on tail}
Next, we can apply the scheme to the bump-on-tail instability test case 
for which the initial condition writes 
$$
f_0(x, v) = \tilde{f}(v)[1+\alpha\cos(k x)], 
$$ 
with 
$$
\tilde{f}(v)=n_p \exp(-v^2/2) + n_b \exp\left(-\frac{|v-u|^2}{2v^2_t}\right)
$$
on the interval $[0, 20\pi]$, with periodic conditions in space. 
The initial condition $f_0$ is a Maxwellian distribution function which has a bump on 
the Maxwell distribution tail; the parameters of this bump are the following 
$$
n_p=\frac{9}{10(2\pi)^{1/2}}, n_b=\frac{2}{10(2\pi)^{1/2}}, u=4.5, v_t=0.5, 
$$
whereas the numerical parameters are $N_x=128, N_v=128, v_{max}=9, \Delta t=0.5$. 
The Runge-Kutta $4$ algorithm is used to compute the characteristics.  

We are interested in the time evolution of the spatially integrated 
distribution function 
$$
F(t, v) = \int_0^{20\pi} f(t, x, v) dx, 
$$
and in the time history of the electric energy $1/2 \|E(t)\|^2_{L^2}$.   
For this latter diagnostic, we expect oscillatory behavior of 
period equal to $1.05$; moreover, since an instability will be declared, 
the electric energy has to increase up to saturation at $t\approx 20.95$ 
and to converge for large times to $36 \%$ of its highest value 
(see \cite{nakamura, shoucri}).  

\begin{figure}
\psfrag{time}{$t \, \omega_p^{-1}$}
\psfrag{EE}{Electric energy}
\begin{tabular}{cc}
\includegraphics[width=7cm,height=7.cm,angle=-0]{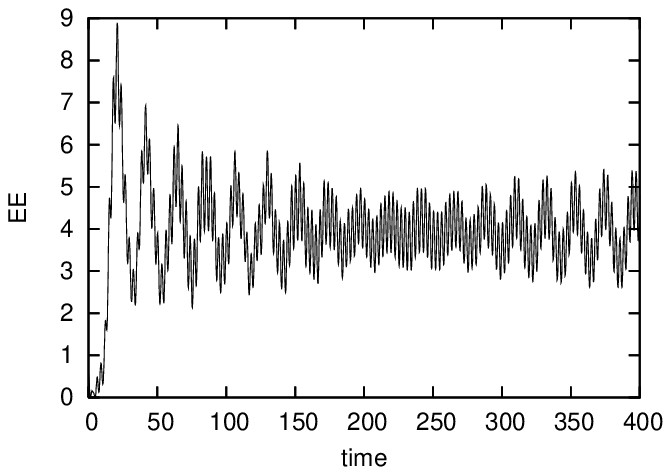} &
\includegraphics[width=7cm,height=7.cm,angle=-0]{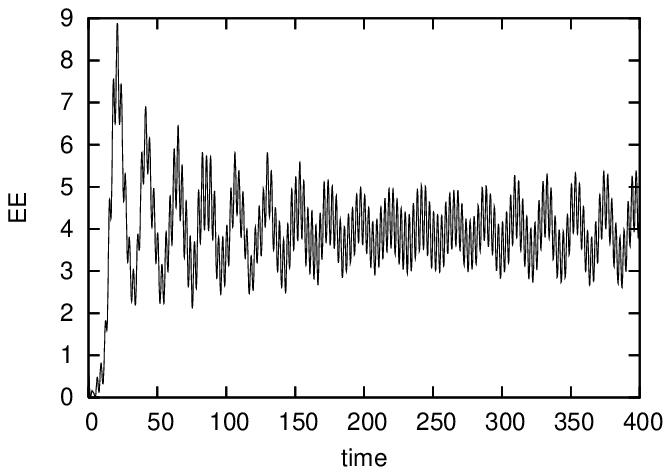}
\end{tabular}
\caption{Bump on tail instability: time evolution of the electric energy for FSL (left) and BSL (right). }
\label{bot1}
\end{figure}

\begin{figure}
\psfrag{EE}{Electric energy}
\psfrag{time}{$t \, \omega_p^{-1}$}
\begin{tabular}{cc}
\includegraphics[width=7cm,height=7.cm,angle=-0]{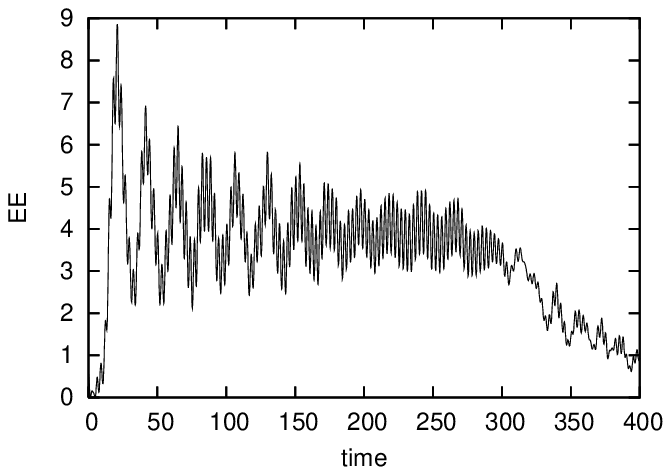} &
\includegraphics[width=7cm,height=7.cm,angle=-0]{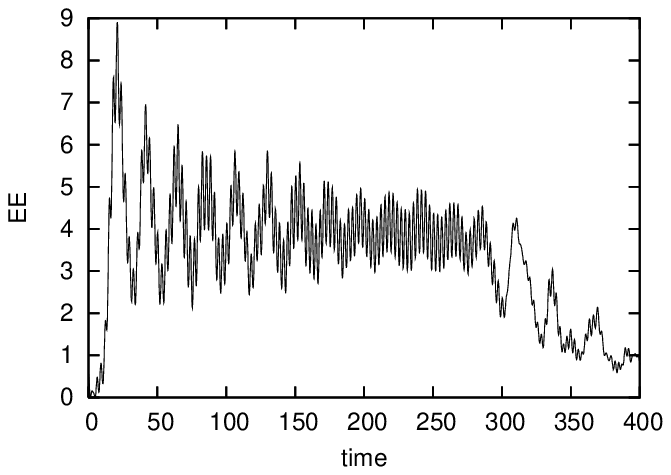}
\end{tabular}
\caption{Bump on tail instability: time evolution of the electric energy for FSL RK2 (left) and FSL Verlet (right). }
\label{bot11}
\end{figure}

On Figures \ref{bot1} and \ref{bot11}, we plot the electric energy as a function of time. 
We can observe that oscillations appear, the period of which can be evaluated 
to $1.$; then the maximum value is reached at $t \approx 21$ and the 
corresponding amplitude is about $9$, which is in very good agreement 
with the results presented in \cite{nakamura}. Then the amplitude of the 
electric energy decreases and presents a slower oscillation due to the 
particle trapping. Finally, it converges to an amplitude of about $2.8$ 
which is very close to the predicted value. 
However, for very large times 
(at $t=250\,\omega_p^{-1}$), FSL using lower time integrator algorithms 
(Runge-Kutta $2$ or Verlet algorithms) leads to bad results 
(see Fig. \ref{bot11}). A very precise computation of the characteristics 
is required for this test and the  use of Runge-Kutta $4$ is crucial. 
The results obtained by BSL (see Fig. \ref{bot1}) are very close to those 
obtained by FSL using Runge-Kutta $4$. Moreover, if $\Delta t$ increases 
(up to $\Delta t\approx 0.75$), BSL presents some unstable results whereas FSL 
remains stable up to $\Delta t \approx 1$. The use of high order time integrators  
leads to an increase of accuracy but also to more stable results as $\Delta t$ 
increases. It has been remarked that linear interpolation of the electric field 
is not sufficient to obtain accurate results with FSL-Runge-Kutta $4$ and cubic spline 
are used to that purpose.  

\begin{figure}
\psfrag{fv}{Integ. d.f.}
\psfrag{v}{velocity}
\begin{tabular}{cc}
\includegraphics[width=7cm,height=7.cm,angle=-0]{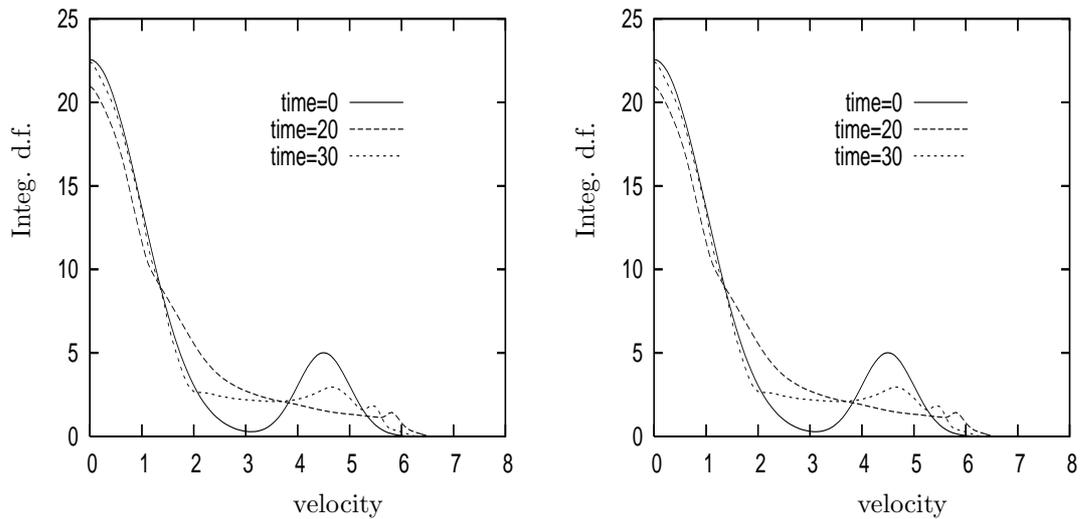} &
\includegraphics[width=7cm,height=7.cm,angle=-0]{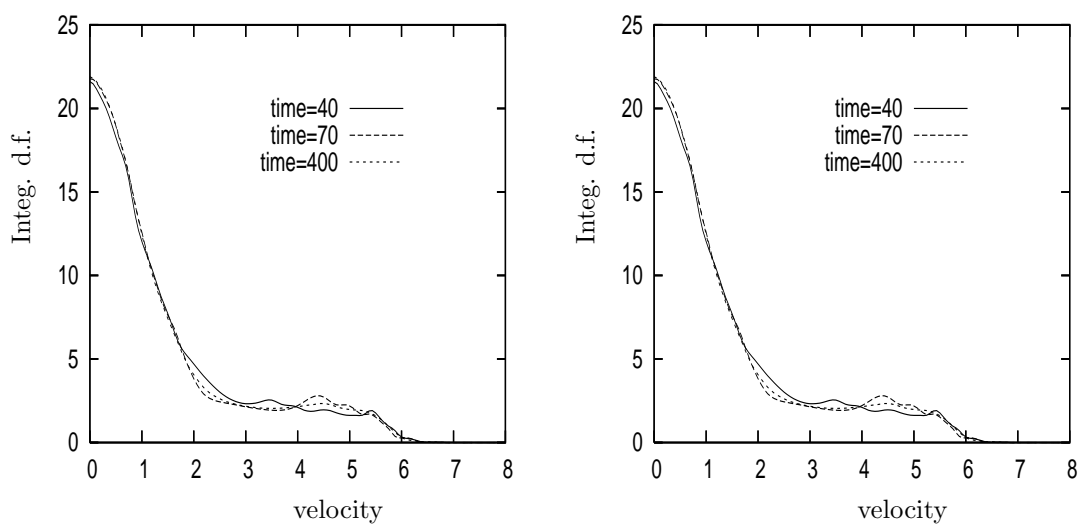}
\end{tabular}
\caption{Bump on tail instability: time development of the spatially integrated distribution function for FSL (left) and BSL (right). }
\label{bot2}
\end{figure}

\begin{figure}
\psfrag{fv}{Integ. d.f.}
\psfrag{v}{velocity}
\begin{tabular}{cc}
\includegraphics[width=7cm,height=7.cm,angle=-0]{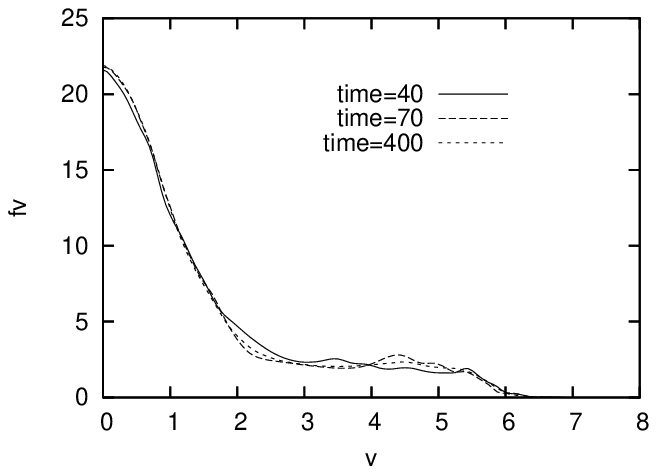} &
\includegraphics[width=7cm,height=7.cm,angle=-0]{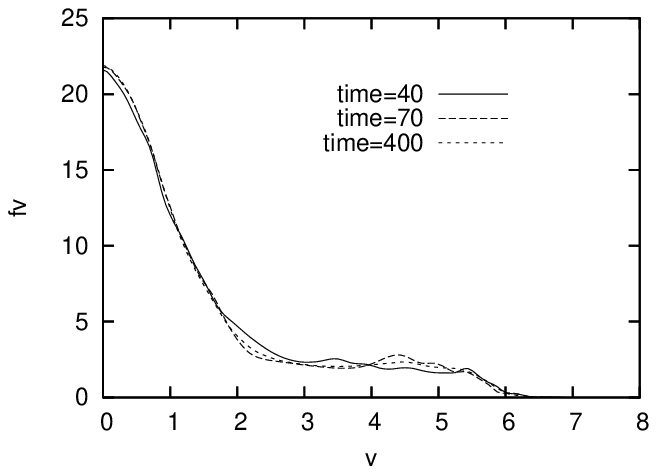}
\end{tabular}
\caption{Bump on tail instability: time development of the spatially integrated distribution function for FSL (left) and BSL (right). }
\label{bot3}
\end{figure}

Fig. \ref{bot2} and Fig. \ref{bot3}  shows the time development of the 
spatially integrated distribution function for FSL and BSL. 
We observe that very fast, the bump begins to be merged by the 
Maxwellian and a plateau is then formed at $t\approx 30-40 \, \omega_p^{-1}$. 

\subsection{Guiding-center case}
\paragraph{Kelvin-Helmoltz instability}

In order to validate our guiding-center code, we used two test cases introduced in \cite{shoucri}  and \cite{ghizzo}. 

\subsubsection{First test case}

The corresponding initial condition is 
$$
\rho(x,y,t=0) = \rho^0(y) + \epsilon \rho^1(y)
$$
coupled with Poisson's equation:
$$
\phi = \phi^0(y) + \epsilon \phi^1(y) \cos(kx)
$$
The instability is created choosing an appropriate $\rho^1$ which will perturb the solution around the equilibrium one $(\rho^0,\phi^0)$. Using the the work of Shoucri, we will take:
$$
\rho(x,y,t=0) = \sin(y) + 0.015\sin(\frac{y}{2})\cos(kx)
$$
where $k=\frac{2 \pi}{L_x}$ and $L_x$ the length of the domain in the $x$-direction 

The numerical parameters are:
$$
N_x = N_y = 128, \Delta t = 0.5. 
$$
The domain size has an impact on the solution. The interval $[0,2 \pi]$ 
will be used on the $y$-direction, and respectively $L_x=7$ and $L_x=10$ . 
This leads to real different configurations:

With $L_x=7$, Shoucri proved that the stable case should be dealt with. 
That is what was observed with this code.

\begin{figure}
\begin{tabular}{cc}
\includegraphics[width=7cm,height=7.cm,angle=-0]{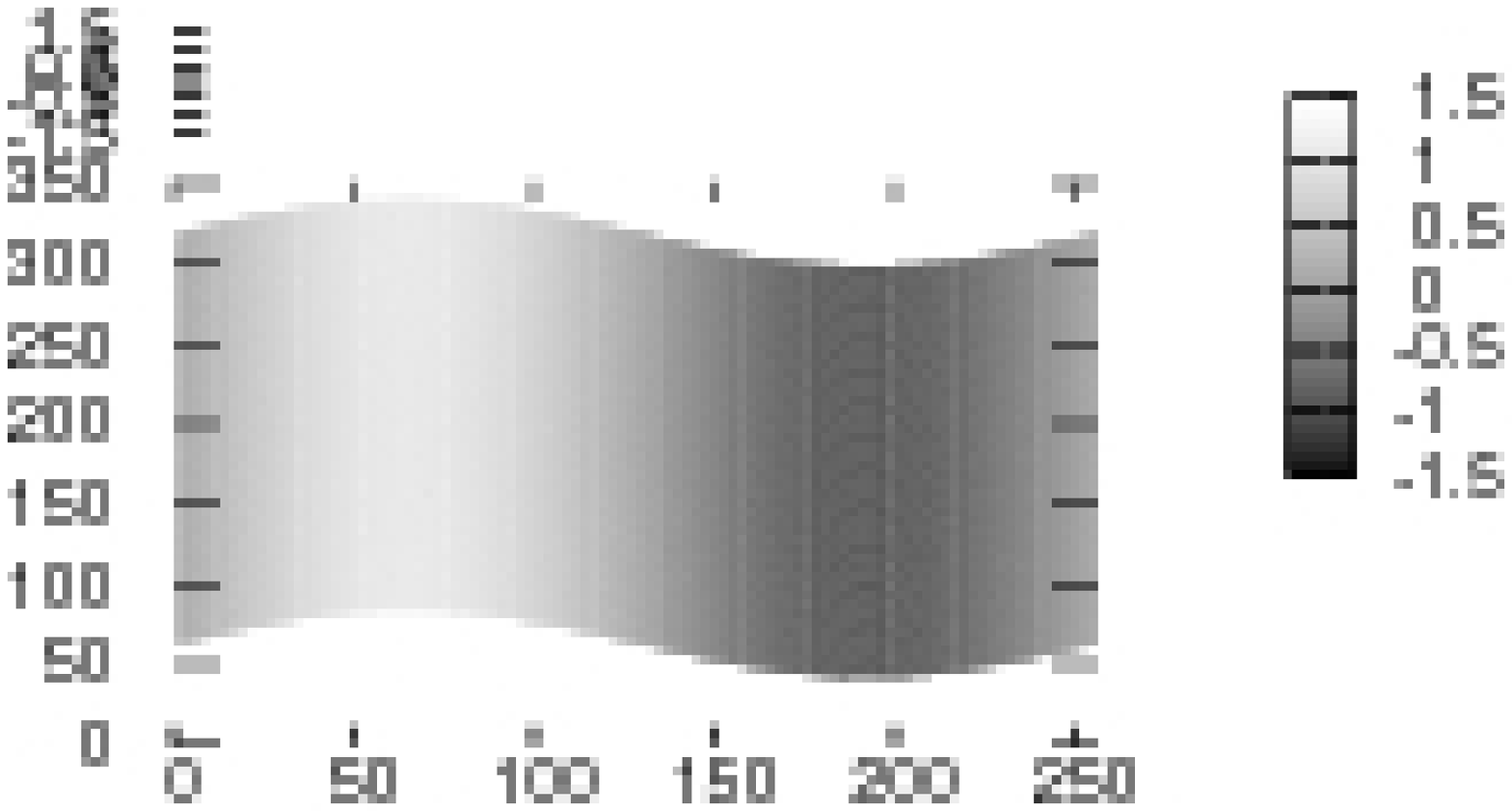} &
\includegraphics[width=7cm,height=7.cm,angle=-0]{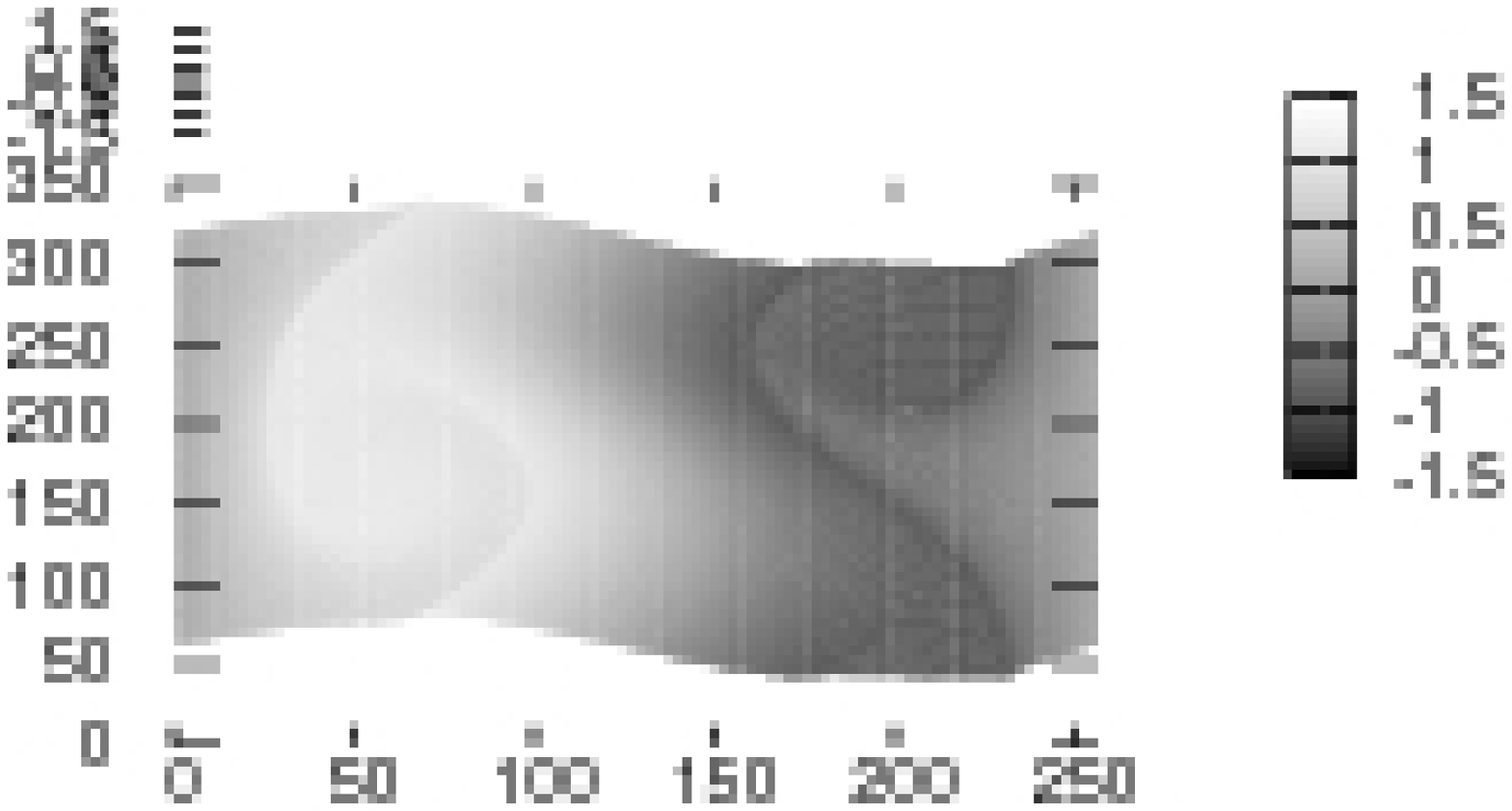}
\end{tabular}
\caption{Kelvin Helmholtz instability 1: distribution function at time $t=0, 30\, \omega_p^{-1}$}
\label{KH1_1}
\end{figure}

\begin{figure}
\begin{tabular}{cc}
\includegraphics[width=7cm,height=7.cm,angle=-0]{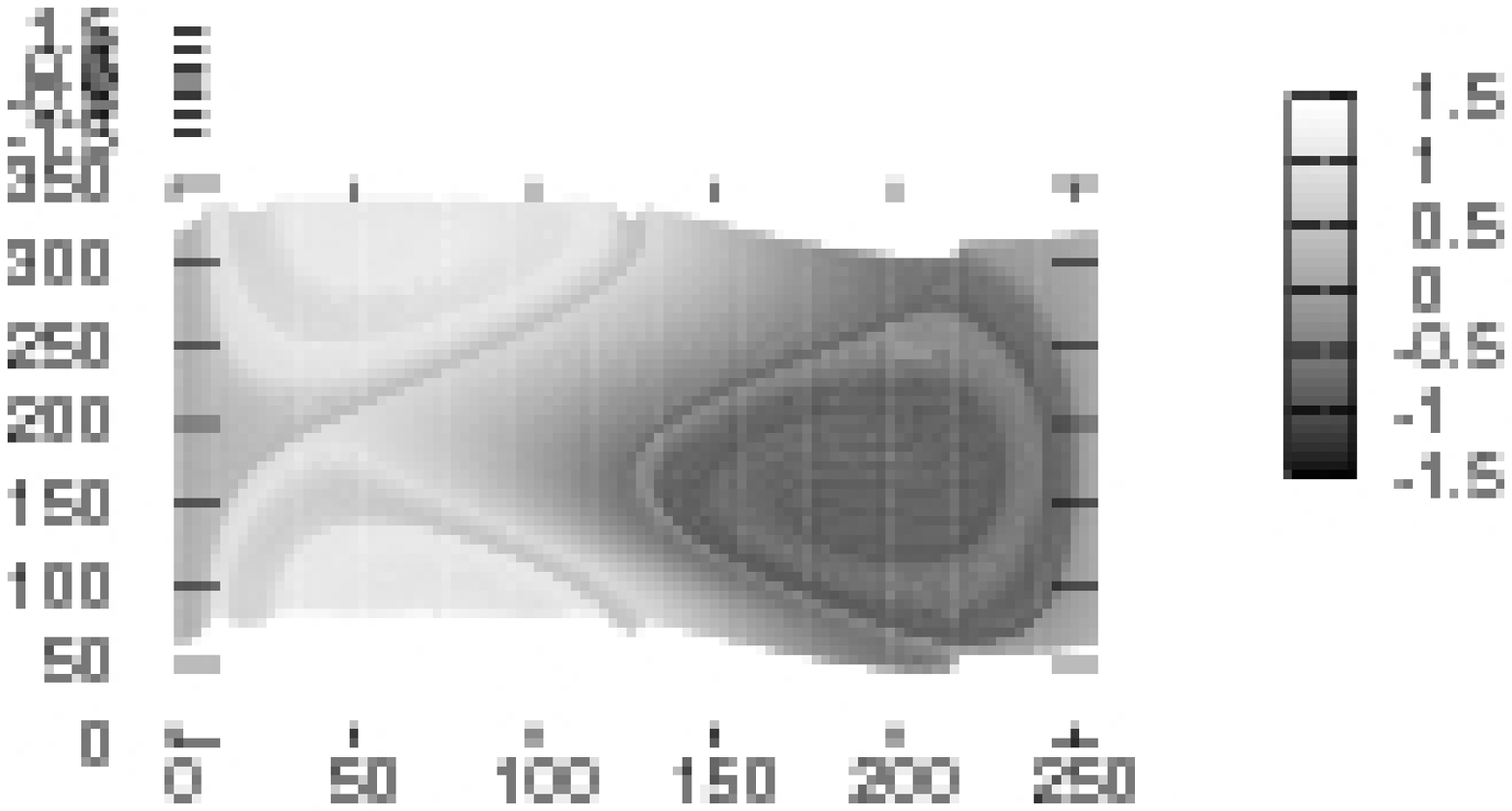} &
\includegraphics[width=7cm,height=7.cm,angle=-0]{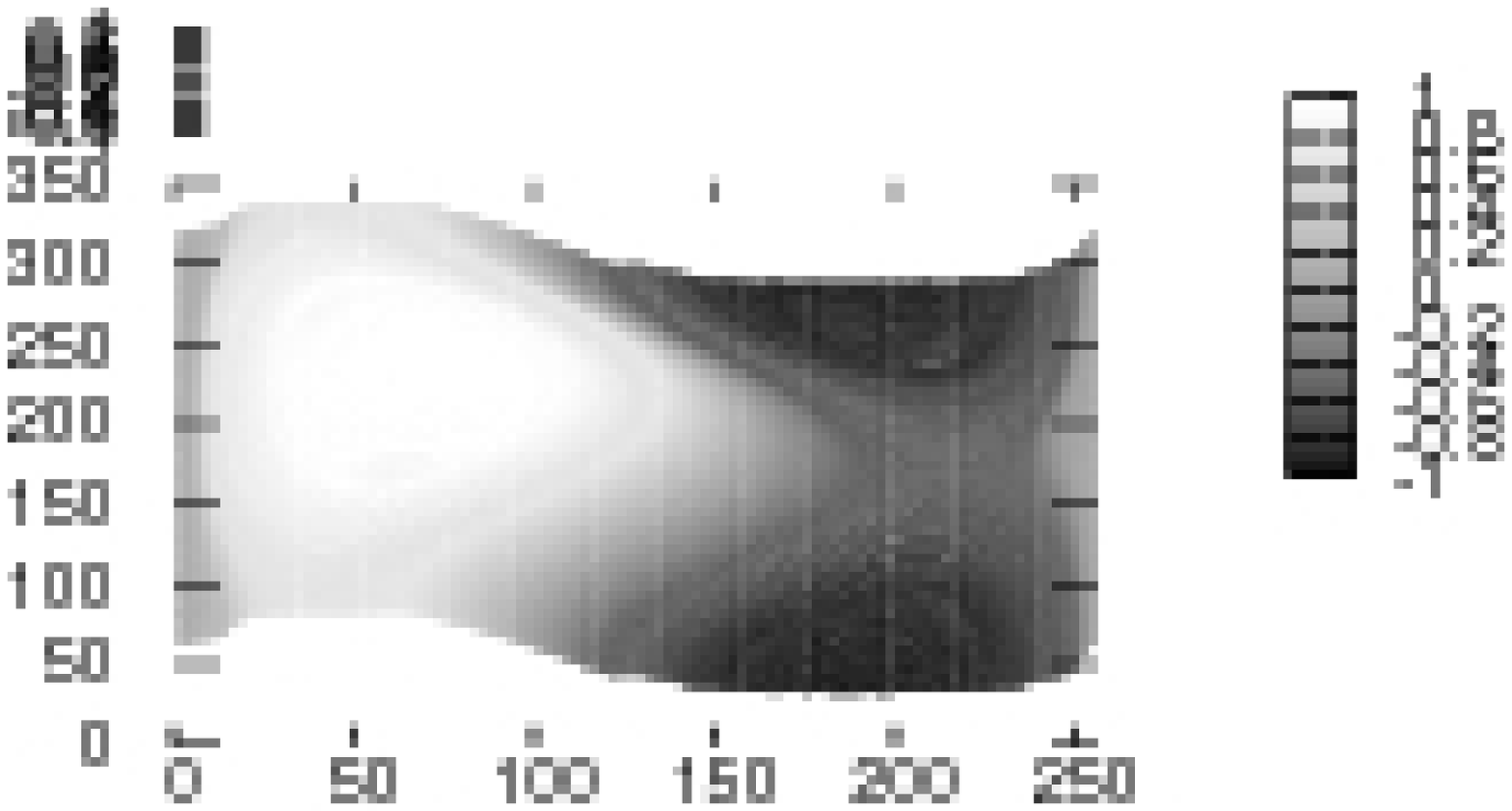}
\end{tabular}
\caption{Kelvin Helmholtz instability 1: distribution function at time $t=50, 500\, \omega_p^{-1}$}
\label{KH1_2}
\end{figure}

With $L_x=10$, the unstable case is faced. The results prove it on 
figure Fig. \ref{KH1_1} and \ref{KH1_2}. 

\begin{figure}
\psfrag{time}{$t \, \omega_p^{-1}$}
\psfrag{ord}{$\|E\|_{L^2} $}
\begin{tabular}{cc}
\includegraphics[width=7cm,height=7.cm,angle=-0]{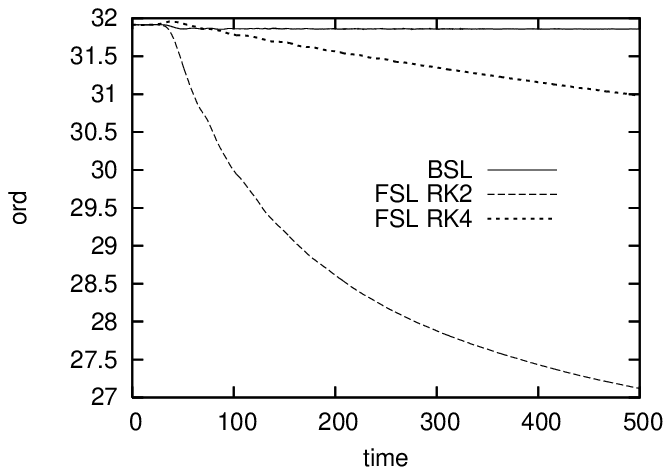} &
\psfrag{ord}{$\|\rho\|_{L^2} $}
\includegraphics[width=7cm,height=7.cm,angle=-0]{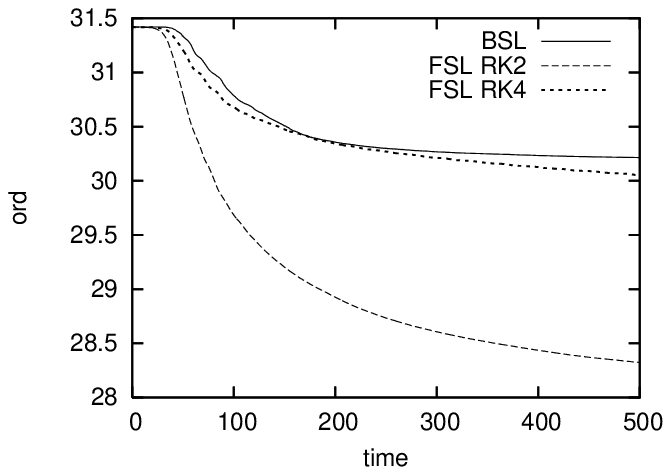}
\end{tabular}
\caption{Kelvin Helmholtz instability 1: time history of $L^2$ norms of $E$ (left) and of $\rho$ (right). Comparison between FSL and BSL. }
\label{KH1_4}
\end{figure}

\begin{figure}
\psfrag{delta_t=0.25}{$\Delta t=0.25\;$}
\psfrag{delta_t=0.5}{$\Delta t=0.5\;$}
\psfrag{delta_t=0.75}{$\Delta t=0.75\;$}
\psfrag{delta_t=1}{$\Delta t=1\;$}
\psfrag{ord}{$\|E\|_{L^2} $}
\begin{tabular}{cc}
\includegraphics[width=7cm,height=7.cm,angle=-0]{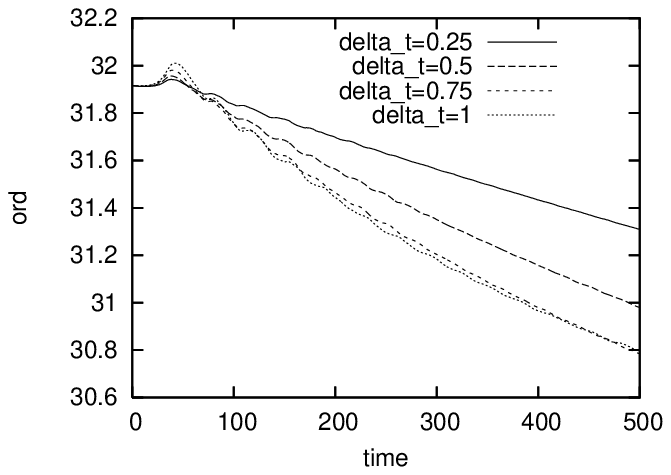} &
\psfrag{ord}{$\|\rho\|_{L^2} $}
\includegraphics[width=7cm,height=7.cm,angle=-0]{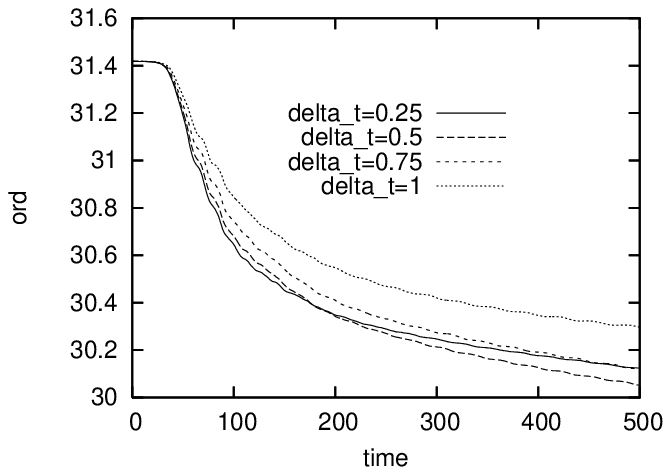}
\end{tabular}
\caption{Kelvin Helmholtz instability 1: time history of $L^2$ norms of $E$ (left) and of $\rho$ (right). Comparison of the results for RK4 for different values of $\Delta t$. }
\label{KH1_5}
\end{figure}

For this test case, the evolution of the energy $\int E^2 dx dy$ and 
enstrophy $\int \rho^2 dx dy$ will also be plotted on Fig. \ref{KH1_4}. 
These should be theoretically invariants of the system. Like for 
other semi-Lagrangian methods, the energy lowers during the first 
phase, which is the smoothing one, where micro-structures can not be 
solved properly. Nevertheless, the energy is well conserved.
Moreover, on Fig. \ref{KH1_4}, FSL using second and fourth order 
Runge-Kutta's methods are compared to the BSL method. 
As observed in the bump-on-tail test, the Runge-Kutta $4$ leads to 
more accurate results and is then very close to BSL. 
However, BSL seems to present slightly better behavior 
compared to FSL-Runge-Kutta $4$. But FSL-Runge-Kutta $4$ enables to simulate 
such complex problems using higher values of $\Delta t$ (see Fig. \ref{KH1_5}).  
We observed for example that the use of $\Delta t=1$ gives rise to very reasonable 
results since the $L^2$ norm of the electric field $E$ decreases of about 
$4 \, \%$. Let us remark that BSL becomes unstable for $\Delta t \geq 0.7$.

\section{Conclusion and perspectives}
In this paper, we introduced the forward semi-Lagrangian method 
for Vlasov equations. The method has been tested on two 
different models, the one-dimensional Vlasov-Poisson one, and 
the guiding-center one. Different test cases have been simulated, 
and they are quite satisfying. 
The results are in some cases a bit less accurate, with respect to the conservation 
of invariants, than with the classical BSL
method, but enables the use of very large time steps without being unstable 
and recovering all the expected aims. No iterative methods anymore, 
and high order time schemes can be use in a straightforward manner. 
The next step will be to test 
the method with the Vlasov-Maxwell model, in which we will try to 
solve properly the charge conservation problem, which is the ultimate goal. 
We will try to use PIC results about that conservation, for example in \cite{Barthelme}.
We will also try to prove theoretically the convergence of this method.

\section{Appendix I: Linearized Vlasov Poisson and Landau damping}

In classical plasma physics textbooks, only the dispersion relations are
computed for the linearized Vlasov-Poisson equation. However using
the Fourier and Laplace transforms as for the computation of the dispersion
relation and inverting them, it is straightforward to obtain an exact expression
for each mode of the electric field (and also the distribution function if needed).
Note that each mode corresponds to a zero of the dispersion relation.

The solution of the Landau damping problem is obtained by solving the linearized
Vlasov-Poisson equation with a perturbation around a Maxwellian equilibrium,
which corresponds to the initial condition $f_0(x,v)=(1+\epsilon\cos (kx))/\sqrt{2\pi} e^{-v^2/2}$.
Let us introduce the plasma dispersion function $Z$ of Fried and Conte \cite{fried}
$$
Z(\eta)= \sqrt{\pi} e^{-\eta^2}[i-erfi(\eta)],~~
\mbox{ where }
erfi(\eta) = \frac{2}{\pi} \int_0^{\eta}e^{t^2} dt.
$$
We also have, $Z'(\eta) = -2(\eta Z(\eta)+1).$
Then, denoting by $\widehat{E}(k,t)$ the Fourier transform of $E$ 
and by $\widetilde{E}(k,\omega)$ the Laplace transform of 
$\widehat{E}$, the electric field, solution 
of linearized Vlasov-Poisson satisfies:
$$
\widehat{E}(k,t) = \sum_j Res_{\omega = \omega_j} \widetilde{E}(k,\omega) e^{-i\omega t}
$$
where
$$
\widetilde{E}(k,\omega) = \frac{N(k,\omega)}{D(k,\omega)}
$$
$$
D(k,\omega) = 1- \frac{1}{2k^2} Z'(\frac{\omega}{\sqrt{2}k}), \quad
N(k,\omega) = \frac{i}{2\sqrt{2}k^2}  Z(\frac{\omega}{\sqrt{2}k})
$$
The dispersion relation corresponds to
$$
D(k,\omega) =0. 
$$

For each fixed $k$, this equation has different roots $\omega_j$, 
and to which are associated the residues defining $\widehat{E}(k,t)$ that can be computed with Maple. 
These residues that in fact the values 
\be
\label{rd}
\frac{N(k,\omega_j)}{\frac{\partial D}{\partial \omega}(k,\omega)}.
\ee
Let us denote by $\omega_r = Re(\omega_j)$, $\omega_i = Im(\omega_j)$, 
$r$ will be the amplitude of (\ref{rd}) and $\varphi$ its phase. 

\textbf{Remark:}
For each root $\omega=\omega_r + i \omega_i$, linked to $r e^{i\varphi}$, 
there is another: $-\omega_r + i \omega_i$ linked to $r e^{-i\varphi}$. 
Then keeping only the roots in which $\omega_i$ is the largest, which are the dominating ones after a short time, we get:
$$
\widehat{E}(k,t) \approx r e^{i\varphi } e^{-i(\omega_r + i \omega_i)t} + r e^{-i\varphi } e^{-i(-\omega_r + i \omega_i)t} \\
=2r e^{\omega_i t} \cos(\omega_r t - \varphi)
$$
Taking the inverse Fourier transform, we finally get an analytical expression for the dominating mode of the electric field, which we use to benchmark our numerical solution:
$$
E(x,t) \approx  4 \epsilon r  e^{\omega_i t} \sin(kx) \cos(\omega_r t - \varphi)
$$

\textbf{Remark}
This is not the exact solution, because we have kept only the 
highest Laplace mode. Nevertheless, after about one period in time, 
this is an excellent approximation of E, because the other 
modes decay very fast.

\section{Appendix II: Solution of Poisson in the Guided Center model}

\subsection{ Find $\phi$}

We have to solve:
$$
-\Delta \phi (x,y) = \rho(x,y)
$$
We use a Fourier transform in the x direction. This leads, for $i \in [1,N_x]$:
$$
\frac{\partial^2 \widehat{\phi_i}(y)}{\partial y^2} = \xi^2  \widehat{\phi_i}(y) +  \widehat{\rho_i}(y)
$$
Let us introduce a notation. Thanks to Taylor Young formula, we have:
$$
\delta u_i = \frac{u_{i+1}-2u_i+u_{i-1}}{\Delta y^2} = \left(1+\frac{\Delta y^2}{12} \frac{\partial^2}{\partial y^2}\right)\frac{\partial^2 u_i}{\partial y^2} + {\cal O}(\Delta_y^2). 
$$
Let us apply this to $\widehat{\phi_i}$
$$
\delta \widehat{\phi_i} = \left(1+\frac{\Delta y^2}{12} \frac{\partial^2}{\partial y^2}\right) (\xi^2  \widehat{\phi_i} +  \widehat{\rho_i}) + {\cal O}(\Delta_y^2) = \xi^2  \widehat{\phi_i} +  \widehat{\rho_i} + \frac{\Delta y^2}{12}(\xi^2 \delta \widehat{\phi_i} + \delta  \widehat{\rho_i}) + {\cal O}(\Delta_y^2). 
$$
Now, factorizing all of it, we get: 
$$
\widehat{\phi_{i+1}}\left(1 - \frac{\xi^2 \Delta_y^2}{12}\right) + \widehat{\phi_i}\left(-2+\frac{10 \xi^2 \Delta_y^2}{12}\right) + \widehat{\phi_{i-1}}\left(1 - \frac{\xi^2 \Delta_y^2}{12}\right) = \Delta y^2 (\widehat{\rho_{i+1}} + 10 \widehat{\rho_i} + \widehat{\rho_{i-1}}) + {\cal O}(\Delta_y^4). 
$$
This is nothing but the solution of a linear system 
$$
A \widehat{\phi} = R, 
$$ 
where A is a tridiagonal and symmetric matrix
and $R$ is a modified right hand side which allows 
to achieve a fourth order approximation. 


\subsection{Find E}

To compute the electric field from the electric potential, 
we have to solve $E=-\nabla \phi$. To achieve this task, a 
quadrature formula is used. 

In the $x$ direction, which is the periodic one, a third order Simpson method 
is used
$$
\int_{x_{i-1}}^{x_{i+1}} E(x,y) dx = -\phi(x_{i+1},y) + \phi(x_{i-1},y) \approx  \frac{1}{6} E_{i-1}(y) + \frac{1}{6} E_{i+1}(y) + \frac{2}{3} E_{i}(y)
$$
where the $(\phi_i)_i$ is given by previous step. 
There is no problem with extreme values, since the system is periodic.  
We then find the values of the electric field $E$ solving another 
tridiagonal linear system.

Whereas on the $y$-direction, Dirichlet conditions are imposed at 
the boundary. So, we can use the same strategy within the domain, 
but not on the two boundary points. Since we have a third order 
solution everywhere, we want to have the same order there, 
therefore we cannot be satisfied with a midpoint quadrature rule which 
is of second order. So we will add corrective terms, in order to 
gain one order accuracy. Here is how we do this.
$$
\int_{y_{0}}^{y_{1}} E(x,y) dy = -\phi(x,1) + \phi(x,0) \approx \frac{dy}{2}( E(x,0)+E(x,1)) - \frac{dy^2}{12} (\rho(x,1)-\rho(x,0)) + f_{\phi}
$$ 
where 
$$
f_{\phi} = \mathbb{F}_y^{-1}\left(\frac{\Delta_y^2}{12} \xi^2(\widehat{\phi(\xi,1}) - \widehat{\phi(\xi,0)})\right) \; =
\frac{\Delta_y^2}{12} \left[-\partial_{xx}(\phi(x,1)-\phi(x,0))\right]. 
$$
We want to find the precision of this method, thus, we would like 
to evaluate the following difference which is denoted by $A$:
$$
A =\phi(x,0)-\phi(x,1) - \frac{dy}{2}(E(x,0)+E(x,1)) + \frac{dy^2}{12} (\rho(x,1)-\rho(x,0)) - f_{\phi}. 
$$
Using the Poisson equation and Taylor expansion, we have
\begin{eqnarray*}
A+E(x,0)+E(x,1) &=& -\frac{2}{dy}(\phi(x,1)-\phi(x,0)) + \frac{dy}{6}(\rho(x,1)-\rho(x,0)) + \frac{dy}{6}(\partial_{xx}(\phi(x,1)-\phi(x,0))\nonumber\\ 
 &=& -\frac{2}{dy}(\phi(x,1)-\phi(x,0)) - \frac{dy}{6}(\partial_{yy}(\phi(x,1)-\phi(x,0))\nonumber\\ 
&=& -\frac{2}{dy}(\phi(x,1)-\phi(x,0)) - \frac{dy^2}{6} (\partial_{yyy}\phi(x,\xi_1) + {\cal O}(\Delta_y^3)
\end{eqnarray*}
where $\xi_1 \in [y_0,y_1]$. Thus, we finally obtain
$$
A+E(x,0)+E(x,1)=-\frac{2}{dy}(\phi(x,1)-\phi(x,0)) + \frac{dy^2}{6} \frac{\partial^2}{\partial y^2} E(x,\xi_1)+ {\cal O}(\Delta_y^3)
$$

Moreover, classical quadrature theory gives us the existence of 
$\xi_2 \in [y_0,y_1]$ such as
$$
\int_{y_{0}}^{y_{1}} E(x,y) dy = \frac{dy}{2}( E(x,0)+E(x,1)) - \frac{dy^3}{12}  \frac{\partial^2}{\partial y^2} E(x,\xi_2). 
$$
Replacing $E(x,\xi_1)$ by $E(x,\xi_2)$, which is of first order in 
our computation leads to 
$$
A+E(x,0)+E(x,1)=-\frac{2}{dy}(\phi(x,1)-\phi(x,0)) -\frac{2}{dy} \int_{y_{0}}^{y_{1}} E(x,y) dy + (E(x,0)+E(x,1)) + {\cal O}(\Delta_y^3)
$$
so that $A={\cal O}(\Delta_y^3)$ which is what was expected.

%
%
%
%
%
%


\end{document}